\documentclass[12pt]{article}
\usepackage{graphicx}
\usepackage{graphicx}
\usepackage{xcolor}
\usepackage[outdir=./]{epstopdf}
\usepackage{comment}
\usepackage{caption}
\usepackage{subcaption}
\usepackage{float}
\usepackage{array,scalefnt}
\newcolumntype{L}{>{\centering\arraybackslash}m{10cm}}
\usepackage{fullpage}
\usepackage[margin=1in]{geometry}
\usepackage[all]{nowidow}
\usepackage{authblk}

\usepackage[numbers]{natbib}

\begin{document}
 \title{{Simulation-Optimization of Automated Material Handling Systems in a Healthcare Facility}}
\author[1a]{Amogh Bhosekar}
\affil[1]{\small Department of Industrial Engineering, Clemson University}
\author[1b]{Tu\u{g}\c{c}e I\c{s}{\i}k}

\author[2]{Sandra Ek\c{s}io\u{g}lu}
\affil[2]{Department of Industrial Engineering, University of Arkansas} 

\author[1c]{Kade Gilstrap}

\author[3]{Robert Allen}
\affil[3]{Perioperative Services, Prisma Health, Greenville \authorcr 
Emails: \{abhosek,tisik,kadeg\}@clemson.edu, sandra@uark.edu, robert.allen@prismahealth.org} 

\date{}

\maketitle
\begin{abstract}
{
Automated material handling systems are used in healthcare facilities to optimize material flow, minimize workforce requirements, reduce the risk of contamination, and reduce injuries. This study proposes a framework that integrates data analysis with system simulation and optimization to address the following research questions: \textit{(i)} What are the implications of redesigning a hospital’s material handling system? \textit{(ii) }What are the implications of improving a hospital’s material handling process? This paper develops a case study using data from the Greenville Memorial Hospital (GMH) in South Carolina, USA. The case study is focused on the delivery of surgical cases to operating rooms at GMH via Automated Guided Vehicles (AGVs). The data analysis provides distributions of travel times, AGV utilization, and AGV movement patterns in the current system. The results of data analysis are integrated in a simulation-optimization model that incorporates the size of AGV fleet and the corresponding routes to improve system efficiency, increase AGV utilization, and reduce congestion. To address research question\textit{ (i)}, a redesign of AGV pathways is evaluated to determine whether congestion is reduced. For research question \textit{(ii)}, the implementation of a Kanban system is proposed to improve AGV utilization by controlling the number of AGVs used daily, based on the volume of surgical cases. An extensive sensitivity analysis, simulation-optimization experiments, and a pilot study are conducted and indicate that the proposed Kanban system leads to significant reductions in congestion and travel times and increased utilization of AGVs.
}

\end{abstract}



\section{Introduction}\label{intro} 

\textbf{Background:} Inefficiencies in supply chain operations contribute to increasing healthcare costs in industrialized nations. Researchers estimate that logistics-related expenses account as much as 40\% of the operating budgets in hospitals \citep{Dobson2015}. Like in manufacturing settings, automation in healthcare has improved the efficiency of material handling systems and reduced the cost of material flow. With lower operational costs in hospitals, the cost of healthcare for patients can be reduced.

Achieving better patient care is the main goal of any healthcare system. While improving hospitals’ supply chain activities, like material handling, is important, not much research has been conducted in this area \cite{Dobson2015}. Previous work focuses on integrating patient care with supply chain activities to improve the efficiency of healthcare systems while providing the necessary attention to patient needs \cite{Meijboom2011,DeVries2011}. For instance, delivering food, medication, and clean linen to the patients and removing waste in a timely fashion requires resources and coordination. Most of these activities are repetitive, occur several times a day, and use a large portion of a healthcare provider’s labor hours. Using AGVs to handle these activities can allow additional employees to be reallocated to tasks that directly impact patient care.

\noindent {\bf Review of the Literature:}   Several researchers study material handling systems in manufacturing settings. \textcolor{black} {Work by \cite{katevas2001mobile} provides detailed guidelines for design and implementation of robotic systems in healthcare and discussion of prototypes and/or products.} Some of the literature focuses on determining the fleet size of AGVs to be used in a logistics network, which is a challenging task. Four types of approaches are adapted to determine the number of AGVs: (a) calculus-based approaches; (b) deterministic optimization approaches; (c) stochastic optimization approaches, and (d) simulation-based approaches \cite{Choobineh2012}. Early works in calculus-based models focus on empty and loaded travel times of AGVs \cite{Egbelu1987,Fitzgerald1985}. However, the travel time depends on congestion, which is affected by facility layout, the speed of vehicles, and load size \cite{Vis2006}. Work by \cite{Mahadevan1993} uses analytical models to evaluate the impact that increasing the flexibility of AGV routes has on the number of AGVs needed. Regression models are also developed to determine how many AGVs are required based on several factors, such as the number of work centers and route lengths, as well as the number of intersections \cite{Arifin2000}.

The use of optimization methods to improve material handling systems is common in the literature. Deterministic methods, such as integer programming, multi-objective optimization, and mixed integer programming are used to model system dynamics \cite{Sinriech1992,Maxwell1982,Rajotia1998}. For example, a minimum flow algorithm is developed to determine the minimum number of vehicles required at a container terminal \cite{Vis2001}. Another analytical model is developed to calculate an upper and a lower bound on the number of vehicles required \cite{Ji2010}. Some researchers utilize stochastic optimization approaches, such as analytical queueing models that minimize the number of AGVs needed. The steady-state behavior of closed queueing networks can be used to estimate the required fleet size. The results of these methods can be compared with simulation models for validation purposes \cite{Tanchoco1987,Choobineh2012}. A hierarchical queueing approach can also be used to determine how many vehicles are required \cite{Mantel1995}.

Simulation-based approaches are considered time-consuming and costly \cite{Vis2006,Egbelu1987}; however, they can handle the complexities and randomness present in real systems. Thus, simulation-optimization models have been previously used in complex inventory replenishment problems \cite{Jalali2015}, medical supply chain modeling \cite{Niziolek2012}, and fleet-sizing problems to understand system performance. Models based on the idle and wait times of machines, parts, and AGVs, as well as the number and speed of vehicles, are developed to evaluate performance measures \cite{Gobal1991,Lee1990}. Other studies use two-stage approaches for system simulation and evaluation, or they develop case studies to determine fleet size and evaluate the impacts on indicators, such as queue sizes, occupation numbers, and service times \cite{Gebenni2008,Yifei2010}. \textcolor{black}{ In healthcare, simulation models are used to evaluate supply chain performance of AGV-based material handling systems versus manual delivery systems \cite{rossetti2001multi,rossetti1998mobile,chikul2017technology}. These works show that AGV-based material handling systems are economically viable and achieve significant performance gains. Furthermore, previous studies address how a fleet of robots can meet the delivery requirements in hospitals.}

\noindent \textcolor{black}{{\bf Research Questions and Contributions:} The research framework proposed here addresses the following \emph{research questions}: \textit{(i)} What are the implications of \emph{redesigning} the material handling system in a hospital, in terms of efficiency, costs, safety, and ease of implementation?\textit{(ii)} What are the implications of \emph{improving} the material handling process in a hospital, in terms of efficiency, costs, safety, and ease of implementation? The proposed framework uses a simulation-optimization model to identify the number of AGVs required daily and their corresponding routes. This framework is demonstrated through a case study in a hospital. In particular, data collected at the hospital has been used to develop a data-driven, discrete event simulation model that captures the traffic flow of AGV movements based on the business rules that govern these movements.}


\textcolor{black}{This paper makes the following \emph{contributions} to the existing body of literature: \textit{(i)} The proposed research framework presented here enables hospitals to identify what factors impact the performance of the material handling system and to develop solutions that improve its efficiency.
Prior works point to the cost savings and benefits of using simulation to model AGV movements. However, based on our review of the literature, only a few papers discuss the use of AGVs in hospitals \cite{Pedan2017}. \textit{(ii)} The research papers cited here treat the vehicle fleet sizing problem as a tactical issue to be addressed at the design stage, but the problem becomes operational when the focus is on selecting the required number of vehicles from a pool of vehicles on a day-to-day basis. The research presented here addresses the operational level issues associated with fleet size selection that impact AGV movements.}

 
\section{Description of the System} 

The research presented in this paper was conducted in collaboration with GMH, one of the seven campuses of Prisma Health in \textcolor{black}{South Carolina, USA}. GMH provides general inpatient services and specialized treatments for heart diseases and cancer. The hospital also houses the Family Birthplace, the Children’s Hospital, and the Children’s Emergency Center.

\begin{figure}[H]
	\centering
	{\includegraphics[width=\textwidth]{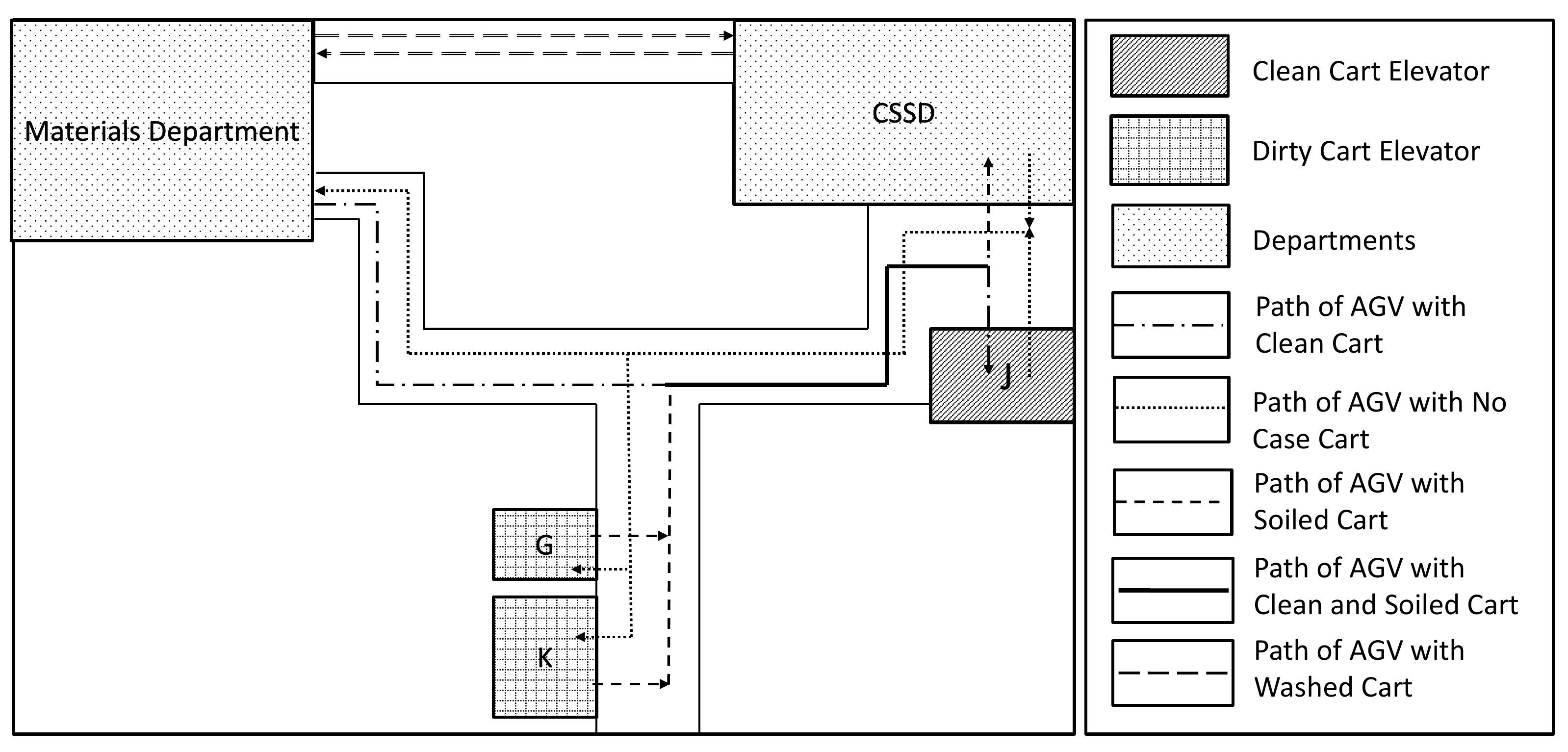}}
	\caption{Map - Mezzanine Floor.}
	\label{fig:Model1}
\end{figure}

This research focuses on material handling activities that support surgical processes at GMH. The research team collaborated with the Perioperative Services Department (PSD) which oversees these processes. The PSD consists of three divisions: the Materials Division (MD) and the Central Sterile Storage Division (CSSD), both of which are located on the mezzanine floor (see Figure \ref{fig:Model1}), and the Operating Room Division (ORD) located on the second floor. The second floor houses 32 operating rooms (ORs) divided into three separate cores. The cores are grouped based on the medical specialties they serve, such as orthopedic treatment, cardiovascular treatment, and neurological treatment. The instruments used in a surgery are stored in the corresponding core.

The type of surgery determines the materials needed, including soft goods and implants, and the surgical instruments used. The PSD is responsible for loading materials to a clean case cart; delivering the case cart from the MD to the OR; loading implants and instruments to clean case carts at the cores; delivering the clean case carts from the cores to the OR; returning soiled instruments, which have already been used, from the OR to the CSSD; and cleaning the soiled instruments at the CSSD. Each case cart is dedicated to a particular surgical case; thus, it contains every material requested by the surgeon. AGVs manufactured by FMC-Technology are used at GMH to move clean and soiled case carts. Figure \ref{fig:Model1} outlines the paths traversed by AGVs and the location of the departments.

The material handling processes managed by PSD are now described. The process starts with the OR manager providing a detailed schedule of the surgeries planned for the next day. Based on the OR schedule and the doctors’ preferences, a list of instruments and soft goods is generated at the MD. Starting at 3 pm, soft goods are loaded manually into clean case carts. This stage is called the \textit{picking process}. Carts are then manually moved to detents. \textcolor{black}{Detents are platforms or areas equipped with the rails necessary for loading and unloading an AGV.} Next, the request for an AGV is submitted via a centralized AGV control system and an available AGV, closest to the MD, is assigned to the case cart. The movement of a loaded AGV is depicted in Figure \ref{fig:Model1} as \emph{``Path of AGV with Clean Cart.”} This AGV uses elevator J to move the cart to the second floor. The clean case cart is then dropped off at \textcolor{black}{one of} the detents in the case cart storage area (CCSA) located next to elevator J on the second floor. 

An inspection takes place at the CCSA to ensure that the soft goods required are delivered. This inspection begins after every case cart has been delivered to the CCSA. If the clean case carts are not delivered by 7 pm, the hospital incurs overtime because employees have to complete their work before they leave for the day. Additionally, GMH uses AGVs to transport dietary and linen carts, both of which have higher priority than surgical case carts after 6 pm. The movement of these carts begins then, and their deliveries are to be completed the same evening. As a result, AGVs become increasingly unavailable for the movement of surgical case carts after 6 pm. Hence, the delivery of surgical case carts should be complete before other services begin requesting AGVs for transport.

The case carts are stored at the CCSA until the next day, the day of the surgery. The instruments and implants, which are stored at \textcolor{black}{one of} the cores, are added to the case cart. The case cart is then moved manually to the OR. After the surgery, the cart is considered soiled and must be decontaminated. A soiled cart is moved manually to the detents on the second floor, and a request for an AGV is submitted to the centralized AGV control system. The assigned AGV moves the dirty cart to the CSSD. The movement of this AGV is depicted in Figure \ref{fig:Model1} as \emph{``Path of AGV with Soiled Cart.”} Then, the AGV uses elevator G or K to move the cart to the mezzanine floor.  The portion of the path that is \emph{shared} by AGVs with clean and soiled case carts is depicted in Figure \ref{fig:Model1} as  \emph{``Path of AGV with Clean and Soiled Cart"}.

The soiled instruments are washed and sterilized at the CSSD to comply with safety guidelines, and the sterilized instruments are loaded to a clean case cart and moved to the corresponding core for storage. The soiled case carts are washed at the cart washer. Once the cart is clean, an automatic request for an AGV is submitted to pick up the washed cart. The movement of AGVs with cleaned carts is depicted in Figure \ref{fig:Model1} as ``\emph{Path of AGVs with Washed Cart"}. Washed carts are dropped at the MD for the picking process. This cycle of surgical case carts starts and ends at the MD, and it is repeated every business day.

For years, the material handling for perioperative service processes at GMH have not changed. However, over the last few years, the number of patients served by GMH has increased rapidly. In an effort to improve the services provided, additional AGVs and case carts were added to the system \textcolor{black}{without updating the physical infrastructure}. As a result, the staff at GMH noticed that AGVs loaded with case carts often sit on the mezzanine floor waiting for elevator J. AGVs coming from the CSSD and elevator J have higher priority than the AGVs moving toward these locations. Hence, AGVs traveling to these locations wait for elevator J \textcolor{black}{for a long time}. Furthermore, AGVs are not allowed to pass each other; thus, if for some reason an AGV stops, the other AGVs following it will also stop at a safe distance, contributing to the traffic. Congestion leads to a shortage of AGVs at the MD. Sometimes, soiled case carts are stuck in traffic, which creates shortage of washed case carts and clean instruments. These shortages lead to further delays in delivering clean carts. Occasionally, the cart washer and instrument washer at the CSSD remain idle for longer periods of time, contributing to underutilization of the equipment.

In 2017, the staff of the PSD reached out to a team of researchers and asked them to investigate the current material handling system and suggest improvements. The GMH team was interested in learning about how a change in AGVs’ current routes and the location of some GMH departments would impact congestion. Based on the layout of the mezzanine floor, \textcolor{black}{it seems intuitive}  that changing the roles of elevators \textcolor{black}{G and K with J} would lead to less congestion. The research team, using the available data, conducted an extensive data analysis of the material handling processes. Based on the results, the research team chose to investigate the potential impacts of reducing the number of AGVs moving surgical case carts. To this end, two simulation models were developed. Section  \ref{DA}  summarizes the data analysis, and Section  \ref{models} provides details of the models developed. Section \ref{simulation} describes the results from our experiments, and Section \ref{results} summarizes the results of model implementation.

\section{Data Collection and Analysis}\label{DA} 

The data collection plan’s main objective is to understand the system, discover inefficiencies, and support discrete event simulation models. Data was collected from the AGV control system \textcolor{black}{for 50 consecutive days}. This data provided information regarding the movement of AGVs, such as date, time, and location of the pickup; date, time, and location of the drop-off; and type of cart an AGV is carrying. Since the scope of the study is restricted to surgical services, only the data on surgical case cart movement is analyzed. In the event of an AGV breakdown, the AGV is moved away from the path of other AGVs and taken to maintenance area. The data points for those AGVs are removed from the data set as outliers. Additionally, AGVs stop if there is any person or another obstacle, such as another AGV, in their scanning radius. \textcolor{black}{Micro-stoppages due to human traffic are not recorded separately in the data and they cause negligible compared to travel times. For this reason, these micro-stoppages are not modeled explicitly. The details on how micro-stoppages due to AGVs traffic is are modeled is explained in the next section as part of our modeling approach.}

\begin{figure}[H]
	\centering
	{\includegraphics[scale=0.4]{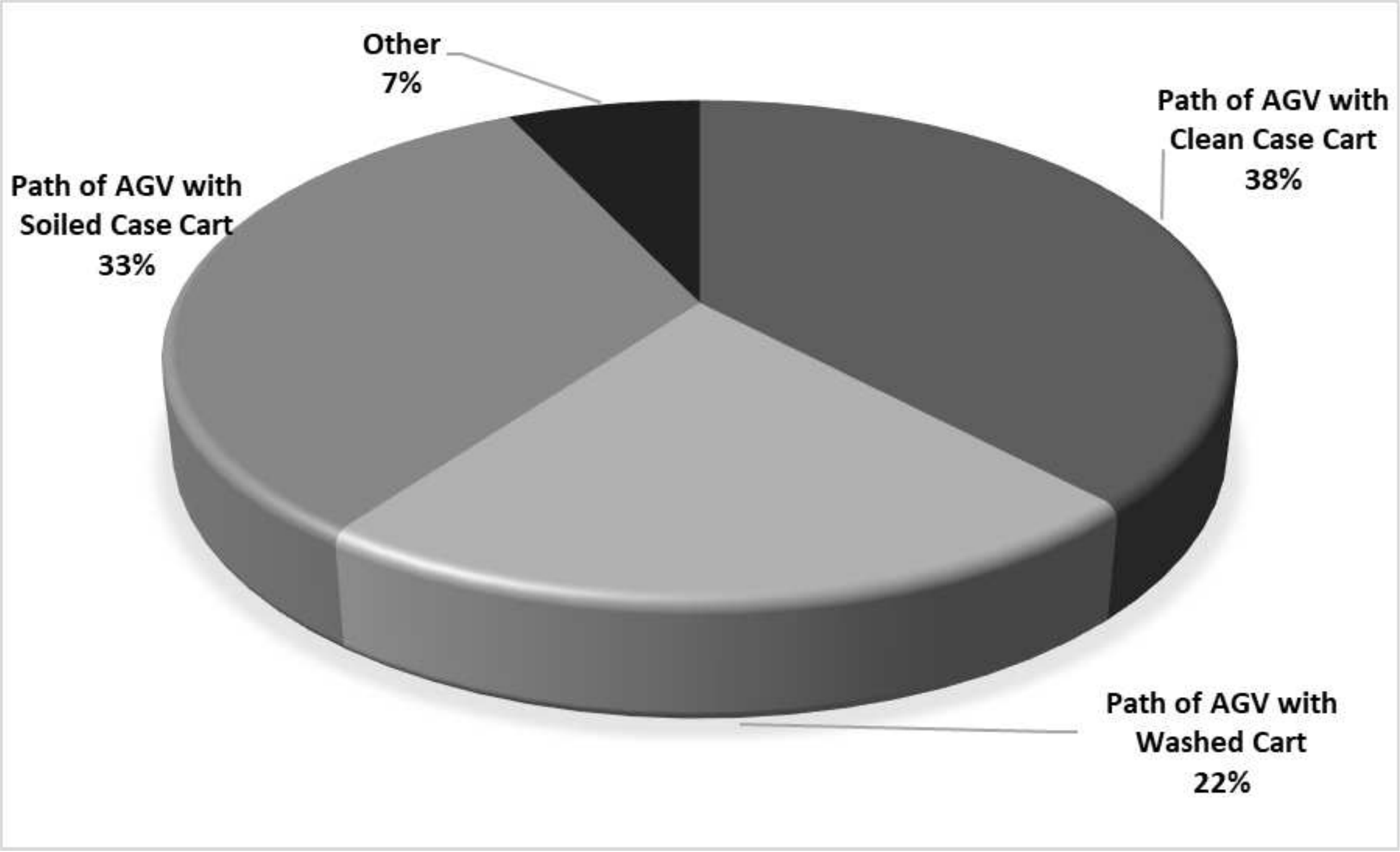}}
	\caption{Number of Trips by Route}    
	\label{fig:PieChart}  
\end{figure}

Figure \ref{fig:PieChart} shows the number of AGV movements along each route. 93\% of the movements are associated with the routes shown in Figure \ref{fig:Model1}, i.e., on the mezzanine floor. Further analysis indicates that 71\% of these movements use the ``\textit{Path of AGV with Clean and Soiled Cart}".  Since the majority of the movements are on these routes, the data analysis focuses on only two routes, ``\textit{Path of AGV with Clean Cart}" and ``\textit{Path of AGV with Soiled Cart}", that contribute most to traffic congestion.

For each cart, the travel time along both routes is calculated based on the data about drop-off and pick-up times. These travel times were grouped in 3 to 4 hour-long intervals. The average and standard deviation of travel times along each path at different times of the day is summarized in Table \ref{Table:ByTime}. These results indicate that travel times are longest along the ``\textit{Path of AGV with Clean Cart}" during 3 pm to 7 pm, \textcolor{black}{when the picking process takes place.}

\begin{table}[ht]
\centering
\caption{AGV Movements by the Time of the Day \label{Table:ByTime}}
{\begin{tabular}{llllll}
\hline
Route & Time & No. of & \multicolumn{2}{c}{Travel Time [Min]} & Coefficient of \\\cline{4-5}
 & Interval & Trips & Average & Std. Dev.  & Variation \\ \hline
 & 12am-3am & 44 & 7.28 & 6.16 & 0.85\\
 & 3am-6am & 53 & 6.87 & 6.16 & 0.9\\
 & 6am-9am & 146 & 5.55 & 3.27 & 0.59\\
2nd Floor Soiled & 9am-12pm & 981 & 4.56 & 4.04 & 0.88 \\
Cart Storage - CSSD & 12pm-3pm & 882 & 5.17 & 2.56 & 0.49\\
 & 3pm-7pm & 753 & \textbf{9.71} & \textbf{8.51} & 0.88 \\
 & 7pm-9pm & 126 & 6.65 & 3.04 & 0.46\\
 & 9pm-12am & 77 & 5.45 & 1.01 & 0.19\\ \hline
 & 12am-3am & 131 & 4.66 & 10.26 & 2.2\\
 & 3am-6am & 227 & 5.96 & 9.55 & 1.6 \\
 & 6am-9am & 112 & 5.27 & 6.17 & 1.17 \\
Materials Department - & 9am-12pm & 80 & 5.8 & 2.4 & 0.41 \\
Case Cart Storage & 12pm-3pm & 101 & 5.33 & 3.69 & 0.69 \\
 & 3pm-7pm & \textbf{1416} & \textbf{8.94} & \textbf{6.49} & 0.73 \\
 & 7pm-9pm & 254 & 5.67 & 4.45 & 0.78\\
 & 9pm-12am & 196 & 4.88 & 5.72 & 1.17\\
 \hline
\end{tabular}}
{}
\end{table}

The results of the data analysis generated the input parameters used in the simulation model. For example, \textcolor{black}{each trip along ``\textit{Path of AGV with Clean Cart}" represents a surgery scheduled for the next day. The total number of surgical cases differs by the day of the week, i.e., Monday to Friday. The number of trips for each day is summed up over each week during this period. This gives seven data points for each day of the week. Because of the limited number of data points, the triangular distribution (TRIA) is used to represent the total number of surgeries scheduled per day. To derive this distribution, the minimum and maximum number of surgeries are determined, and the mode for each day of the week is estimated. The corresponding results are summarized in Table \ref{Table:NrSurgery}. The data for AGV trips on the ``\textit{Path of AGV with Soiled Cart}" is used to estimate the release time of soiled carts from the ORs, since soiled carts are delivered to the CSSD right after a surgery. The data on the total number of soiled carts delivered at the end of every half hour for each day of the week is used to distribute the total number of surgeries over different time intervals within a day. Other input parameters used include the number of AGVs used for surgical case cart movements, the number of case carts available, and the number of cart-washers in the system.}

\begin{table}[ht]
\centering
\caption{Total Number of Surgeries per Day \label{Table:NrSurgery}}
{\begin{tabular}{lc}
\hline
{\bf Day}& {\bf Distribution} \\
\hline
Mon.	& TRIA (60,68,75)\\
Tue.	&TRIA (65,72,76)\\
Wed.	&TRIA (60,65,72)\\
Th.	&TRIA (69,75,80)\\
Fri. &TRIA (55,62,69)\\
 \hline
\end{tabular}}
{}
\end{table}

\textcolor{black}{Note that the existing pathways for movement of AGVs with clean and soiled carts are influenced by safety regulations and the movement of other carts. For example, elevators G and K are used to move soiled surgical instruments, dirty linen, and trash. These elevators continue to the basement to deliver dirty linen and trash. On the other hand, elevator J is solely used for the movement of clean surgical instruments to eliminate any potential contamination. Thus, this elevator serves only the mezzanine and the second floor.}

\section{The Modeling Approach}\label{models}

The simulation model is built using ARENA© simulation software by Rockwell Automation. \textcolor{black}{Tables \ref{Table: inputs} and \ref{Table: Runsetup}  describe the input parameters and run-setup parameters for the simulation model.} A guided path transporter network is constructed to model the movement of AGVs on the mezzanine floor. A guided path transporter functions as a physical entity in the simulation; thus, it is used to model traffic. Transporter-related parameters, such as velocity, acceleration, deceleration, and turning velocity, are obtained from the FMC-Technology AGV handbook. Parameters necessary to model the elevators, such as the time it takes to open and close the door, were obtained from the same handbook. The AGV network consists of intersections and network links. Network links are made up of multiple zones of the same size. AGVs move from one zone to the next along these links. The movement of AGVs is governed by the end control rule, dictating that a transporter releases its current zone at the end of its movement to the next zone. This rule ensures that multiple AGVs can travel on the same network link but not in the same zone. At GMH, the safety distance between AGVs, enforced at all times, is 3 feet. To ensure that AGVs maintain this distance in the model, this study sets a zone length of 3 feet on every network link in the model. \textcolor{black}{When an AGV comes to a halt to maintain sufficient follow-distance or yield to another AGV, it decelerates, momentarily stops, and accelerates again. This study models all these phases of movement, based on the specifications provided in the AGV system handbook. Thus, micro-stoppages due to AGV traffic are modeled accurately in the work presented here.}
\textcolor{teal}{
\begin{table}[ht]
\caption{\textcolor{black}{Input Parameters} \label{Table: inputs}}
{
\begin{tabular}{l|l|l}
\hline
\textbf{Input Parameter} & \textbf{Source} & \textbf{Description} \\ \hline
Number of cases & Surgery Data & Triangular Distribution \\
Soiled case cart release & AGV System Data & Discrete Distribution \\
Number of case cart & AGV System Data & Fixed Capacity: 110 \\
Number of AGVs & AGV System Data & Based on Schedule: \textless{}= 11 \\
Number of cart washers & Employee Survey & Fixed Capacity: 3 \\
Cart washer delay & AGV System Handbook & 20 minutes \\
Clean cart elevator & AGV System Handbook & Fixed Capacity: 2 \\
Soiled cart elevator & AGV System Handbook & Fixed Capacity: 3 \\
Elevator opening/closing delays & AGV System Handbook & 11 seconds \\
AGV network distances & AGV System Floor Maps & Route dependent \\
AGV straight velocity & AGV System Handbook & 200 units distance \\
AGV turning factor & AGV System Handbook & 0.5 \\
Acceleration/Deceleration & AGV System Handbook & 0.98 per second squared \\
Zone control rule & AGV System Handbook & End \\
Case cart loading time & Time Study & TRIA (3,4,5) \\
No. of case cart loading employees & Time Study & Fixed Capacity: 4 \\ \hline
\end{tabular}}%
{}
\end{table}
}

\begin{table}[ht]
\centering
\caption{\textcolor{black}{Run-setup Parameters}
\label{Table: Runsetup}}
{
\begin{tabular}{l|l}
\hline
\textbf{Run-setup parameters} & \textbf{Description} \\ \hline
No. of replications & 30 \\
Base time units & Minutes \\
Warm-up period & 0 \\
Statistics collection & Continuous \\ \hline
\end{tabular}%
}
{}
\end{table}

Destinations such as the CSSD and the MD are modeled using intersections on the network. Each destination is modeled as the last intersection on a path. Network links that connect to these destinations are modeled as bidirectional links. These links have the capacity constraint that at most 1 AGV can be present on the entire link, instead of zones. This constraint ensures that only one AGV can travel to or from any destination on the corresponding link. At every intersection, \textcolor{black}{unless another priority rule applies} as described below, the first-come, first-served rule is followed to determine the right-of-way for the AGVs. If an AGV already has control of an intersection, another AGV must wait to use the intersection until the first AGV leaves it.

\begin{figure}[H]
	\centering
	{\includegraphics[scale=0.6]{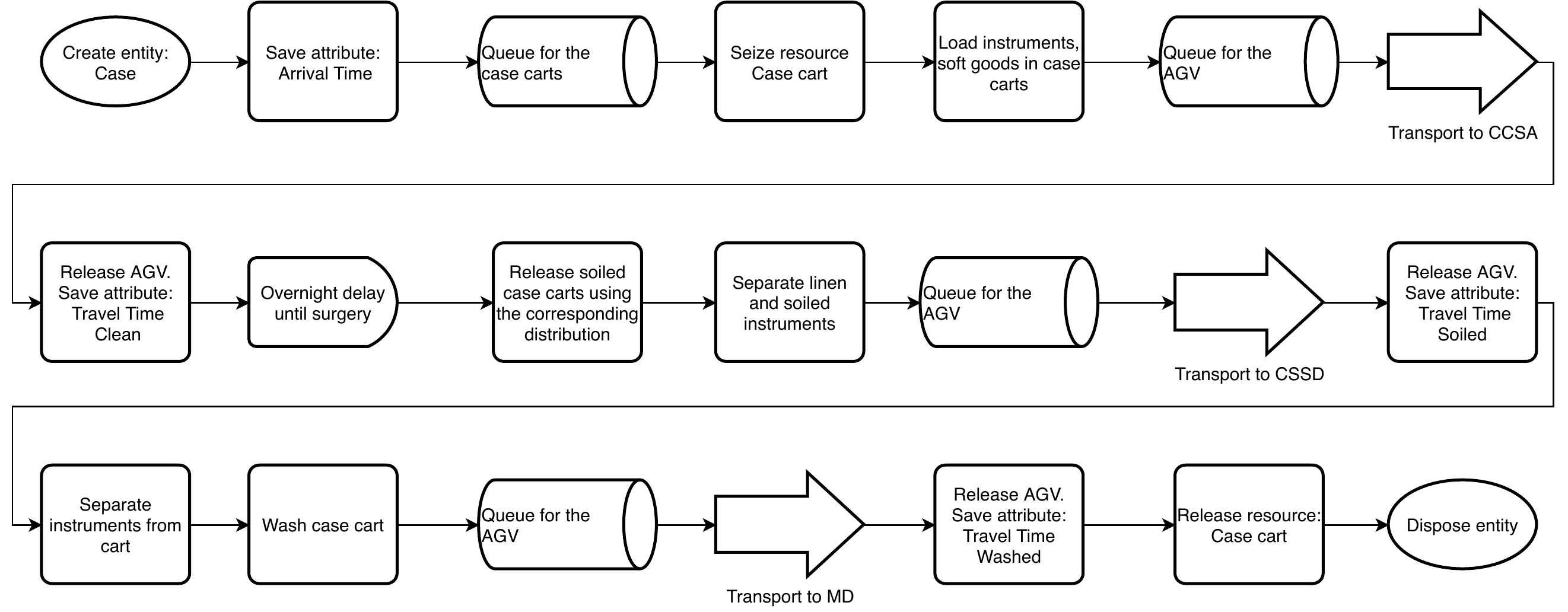}}
	\caption{\textcolor{black}{UML Activity Diagram of Simulation Model} }   
	\label{fig:SimLogic}  
\end{figure}

For the elevators that can accommodate up to 2 AGVs, if there is already an AGV in the elevator, the elevator waits for the next AGV if it is already at the preceding intersection. The length of links in the transporter network are calculated using the GMH floor maps. The detents on the second floor have limited capacity; thus, when they are at full capacity, AGVs are not allowed to enter the elevator because they cannot seize a detent in the second floor. AGVs leaving the departments and/or elevators have higher priority to seize the intersections than the AGVs entering the departments and/or elevators. After completing a task, an AGV is assigned to the next request in the queue. When the queue is empty, AGVs are moved to the parking area. \textcolor{black}{Details of simulation models for the material handling activities are described in Figure \ref{fig:SimLogic} with a unified modeling language (UML) diagram.}

\textcolor{black}{The AGV movements that are only due to surgical cart deliveries are modeled. Other AGV movements are absent in the simulations, thus the AGV traffic due to these deliveries are analyzed in isolation. The models presented here yield accurate results despite these simplifications because \textit{(i)} other AGV activities do not use the same paths as surgical case carts \textit{(ii)} GMH allocates dedicated AGVs to surgical cart movements and AGV availability is not affected by the demand generated for AGVs by other material handling activities.}

\subsection{Movement of Clean Case Carts}

\begin{figure}[H]
	\centering
	{\includegraphics[scale=0.6]{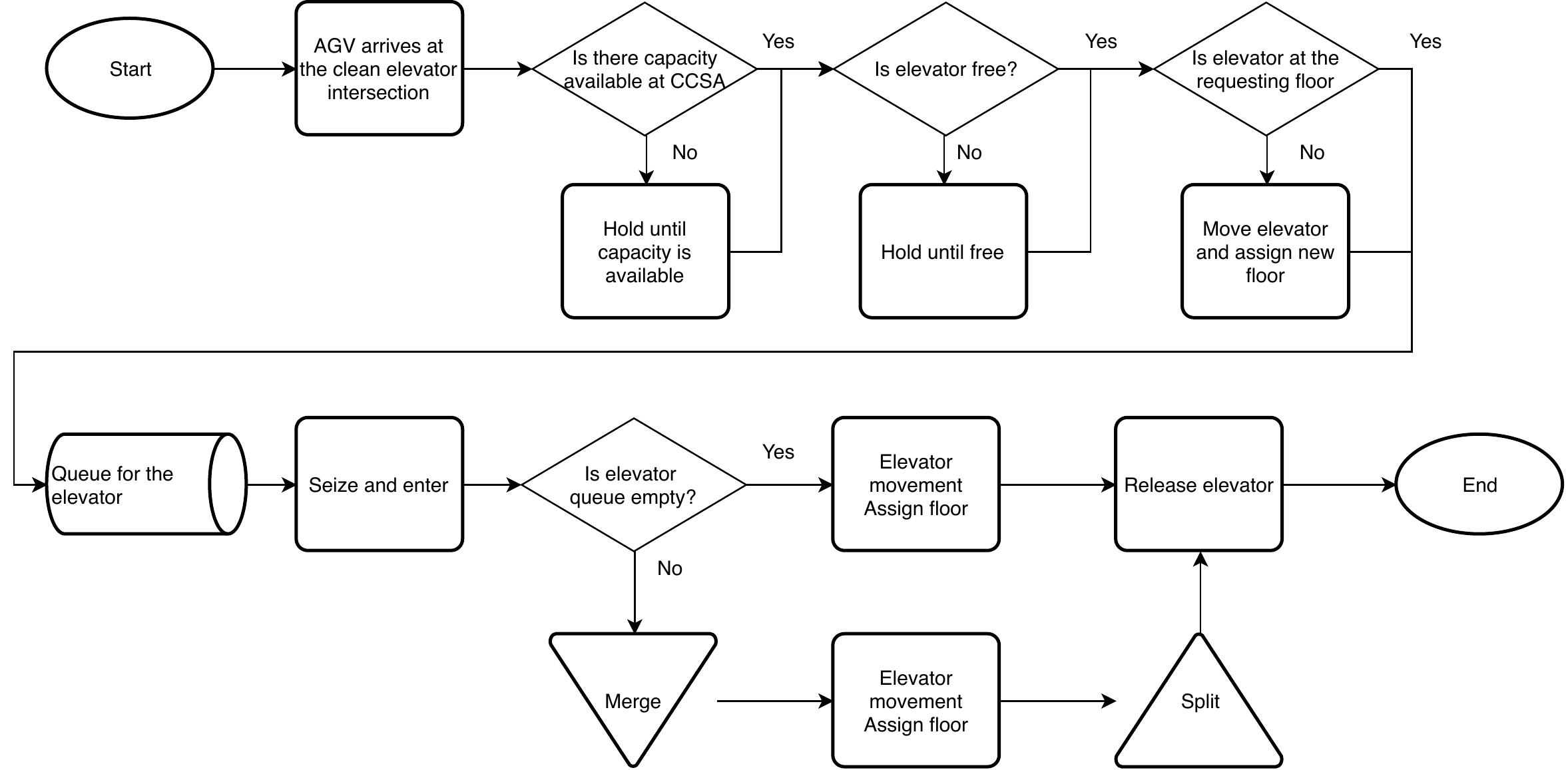}}
	\caption{\textcolor{black}{UML Activity Diagram of Elevator Logic} }
	\label{fig:ElevLogic}  
\end{figure}

Each surgical case is modeled as an entity. Every day, the first entity, whose type is \emph{clean}, is created at 3 pm. Carts and employees loading soft goods to a cart are modeled as resources of the picking process. The distribution of the time required for the picking process is determined to be triangular, TRIA (3,4,5) minutes, based on the data collected via a time study performed. Once the soft goods are loaded into the cart, a request for an AGV is submitted. An available AGV closest to the MD is assigned to the case cart. The AGV travels to the pickup location following the rules governing the AGV’s movement. After picking up a case cart, the AGV travels to elevator J using ``\textit{Path of AGV with Clean Cart}", as shown in Figure \ref{fig:Model1}. A predefined look-ahead stop is modeled before the intersection, ``\emph{Intersection J}", in the corridor between the CSSD and elevator J. At this stop, the availability of intersection J and elevator J, as well as the capacity to accommodate a vehicle at the CCSA, the destination, are checked. If these conditions are satisfied, then the AGV seizes intersection J and advances to its final destination. Otherwise, the AGV is put on hold at the predefined stop until these conditions are satisfied. The clean case cart is dropped at the CCSA. \textcolor{black}{Details of the elevator logic in our simulation models are described in Figure \ref{fig:ElevLogic} with the UML diagram below}. \textcolor{black}{Clean Cart movements are completed by 7 pm every day, and by midnight if there is a need for overtime. Thus, clean cart movements across different days do not interfere.}

\subsection{Movement of Soiled Case Carts}

At 8 am the next day, the carts are released from the CCSA and moved to the ORs to prepare for surgery. The carts are released from surgeries based on the discrete distribution of release times. After they are released from the OR, a new entity type, called \emph{soiled}, is assigned to the case cart. The soiled carts are moved to a location near elevators G and K on the second floor, the \textit{Soiled Cart Storage Area} (SCSA). A request for an AGV is submitted, and an available AGV closest to the SCSA is assigned to the case cart. Soiled case carts are then transported to the CSSD along the ``\textit{Path of AGV with Soiled Cart}" shown in Figure \ref{fig:Model1}. \textcolor{black}{Soiled case cart movements start at 8 am and continue throughout the day based on the surgery schedule. These movements end by 8 am of the next day (see the release time distribution we report in the Appendix.) Thus, soiled case cart movements across different days do not interfere. As a result, there is no accumulation of workload in the system. Given these characteristics, the only initialization effect in the model is due to the lack of dirty cart movement in the very first day of each replication. Observations show that this effect is minimal and does not affect the simulation results.}


\subsection{Movement of Washed Case Carts}

The cart-washer is modeled as a resource with a fixed cycle time of 15 minutes. After an AGV drops a cart at the CSSD, soiled instruments are separated from the cart, and the cart is loaded into the cart washer. After the cart is cleaned, the closest available AGV picks up the cart and transports it to the MD along ``\textit{Path of AGV with Washed Cart}". The carts are subject to an additional 30 minutes of drying time before they are released from the associated surgical case.

\section{Simulation Experiments}\label{simulation}
\noindent \textcolor{black}{This study simulates the following system for 30 days: Each day begins at 8 am and ends at 8 am the next day, corresponding with the actual operating hours of the CSSD, where the clean surgical carts are loaded and the soiled case carts are cleaned.}

\subsection{Validation}

To validate the model presented here, first, a statistical comparison of the current system with simulation model M is conducted based on travel times for clean and soiled case carts. Model M uses 11 AGVs, which is the the same number of AGV's as the hospital \textcolor{black}{currently uses every day} for the movement of case carts. The purpose of this comparison is to ensure that model M accurately reflects the logic and business rules of the current system. \textcolor{black}{The simulation was ran for 30 replications to generate the data necessary for the statistical analysis}. The travel times in the validation runs of model M and the data for clean case cart and soiled case cart movements are tested, using a t-test at a 95\% confidence level, to determine whether they are significantly different. Tables \ref{Table:ValidationCM} and \ref{Table:ValidationDM} summarize the results of the $t$-test. Based on the confidence intervals, it is concluded that the difference between average travel times is not statistically significant, and thus, the simulation model presented here is valid.  These results were also verified by a team from the PSD at GMH.

\begin{table}[h]
\centering
\caption{Model Validation: Clean Case Carts}
\label{Table:ValidationCM}
{\begin{tabular}{ccc}
\hline
\bf{Systems} & \bf{Avg. Travel Time} & \bf{Confidence Interval} \\ \hline
Current System (Data)  & 9.62& (9.36,9.87) \\ 
Model M& 9.67& (9.65,9.68) \\ \hline
\end{tabular}}
\end{table}

\begin{table}[h]
\centering
\caption{Model Validation: Soiled Case Carts
\label{Table:ValidationDM}}
{\begin{tabular}{ccc}
\hline
\bf{Systems} & \bf{Avg. Travel Time} & \bf{Confidence Interval} \\ \hline
Current System (Data) & 6.636& (6.17,6.55) \\ 
Model M & 6.229& (6.217,6.241) \\ \hline
\end{tabular}}
{}
\end{table}

\subsection{Research Question 1: Redesign the System by Swapping the Elevators for Clean And Soiled Case Carts}

\begin{figure}[h]
    \centering
    \includegraphics[width=\textwidth]{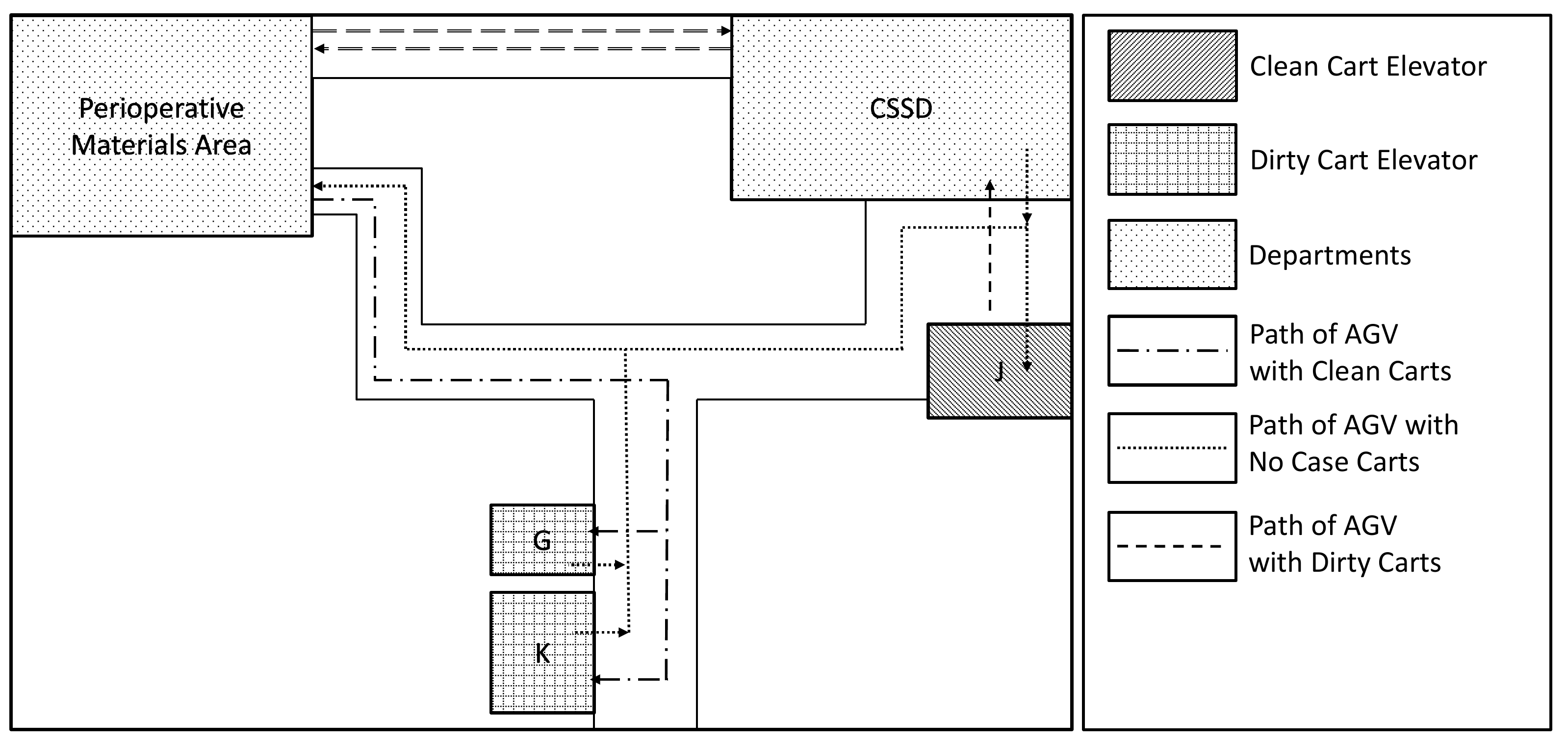}
    \caption{AGV System with Swapped Elevators (S)} \label{fig:Model2} 
        {}
\end{figure}

\textcolor{black}{Swapping the role of elevator J with G and K would eliminate shared paths. Thus, GMH staff believed that swapping elevators would reduce congestion on the mezzanine floor. However, there is a tradeoff between congestion and travel distances since the new design modifies the paths for clean and soiled carts. We have modeled and simulated this design change to evaluate the resulting tradeoff.} Figure \ref{fig:Model2} presents the alternative AGV routes for the clean and soiled case carts when the elevators are swapped. In this case, AGVs that follow the new clean case cart route take elevators G or K and drop the clean case carts in the current system’s SCSA. After a surgery, the soiled case carts are stored at the CCSA. AGVs carrying soiled case carts take elevator J and travel across intersection J toward the sterilization area, the CSSD. Model S is built to capture these changes.

The input data, the total number of AGVs, and the other parameters in models \textcolor{black}{M and S are} are the same. \emph{Task completion time} ($T_c$) is defined as the difference between the time at which the first clean case cart was picked up and the time at which the last clean case cart of the day was dropped off. Models M and S are compared using travel times and task completion times as performance measures, and a sensitivity analysis is conducted to understand the impacts of the proposed changes on these measures.

A set of simulation experiments were run by varying the number of AGVs in each model, from 3 to 11 increasing with increments of 1, and data was collected for travel and task-completion times. The maximum number of AGVs in the experiments was 11 since GMH also currently uses 11 AGVs. The results of the sensitivity analysis are summarized in the next two subsections.

\subsubsection{Results of the Sensitivity Analysis: Clean Cart Movement}

\begin{figure}[H]
    \centering
        \caption{Sensitivity Analysis of Clean Case Carts: Model M vs Model S}
    ~ 
    \begin{subfigure}[b]{0.45\textwidth}
        \includegraphics[width=\textwidth]{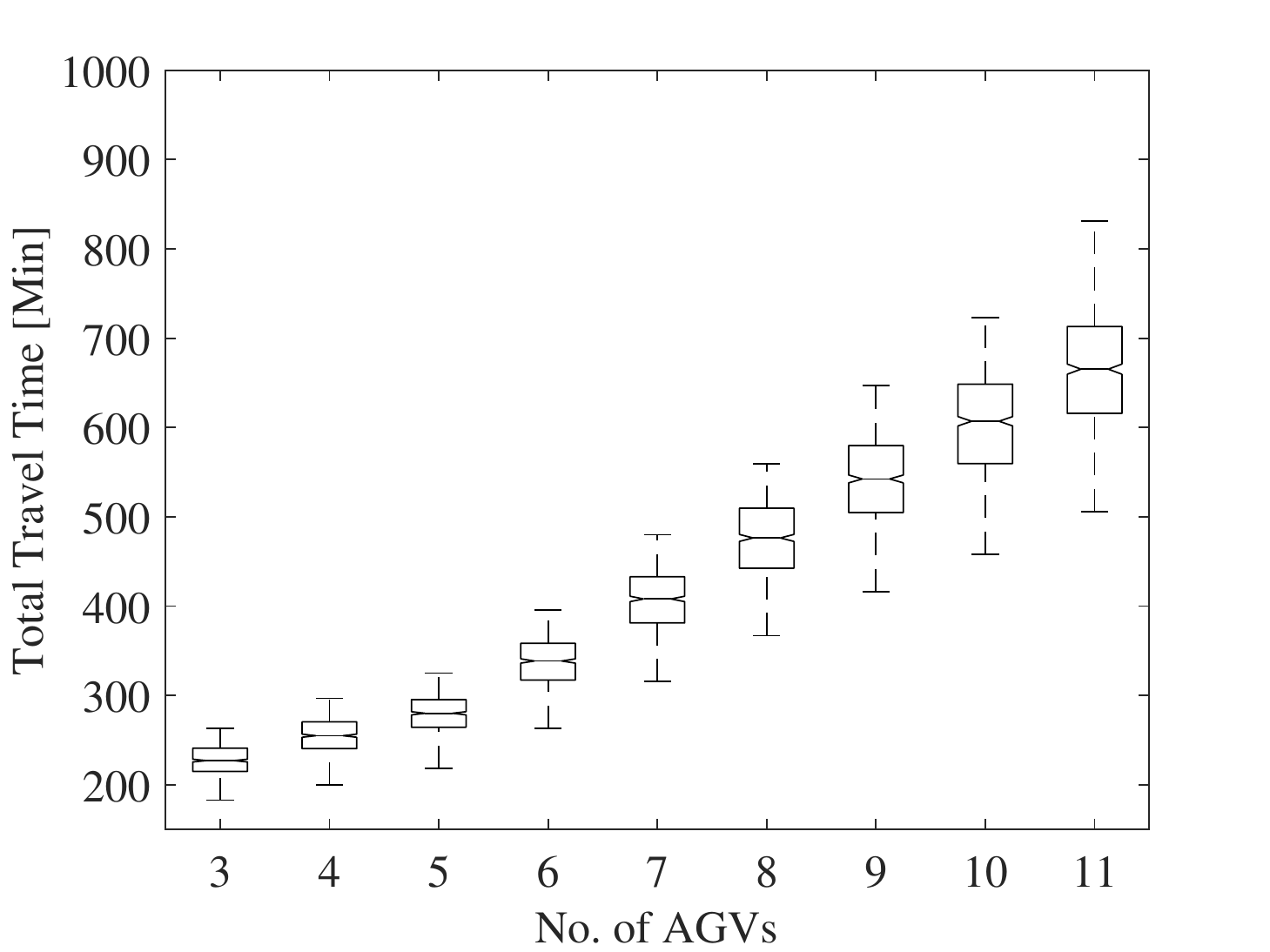}
       \caption{Total Travel Time Model M}
        \label{fig:TTTCM1}
    \end{subfigure}
     \begin{subfigure}[b]{0.45\textwidth}
        \includegraphics[width=\textwidth]{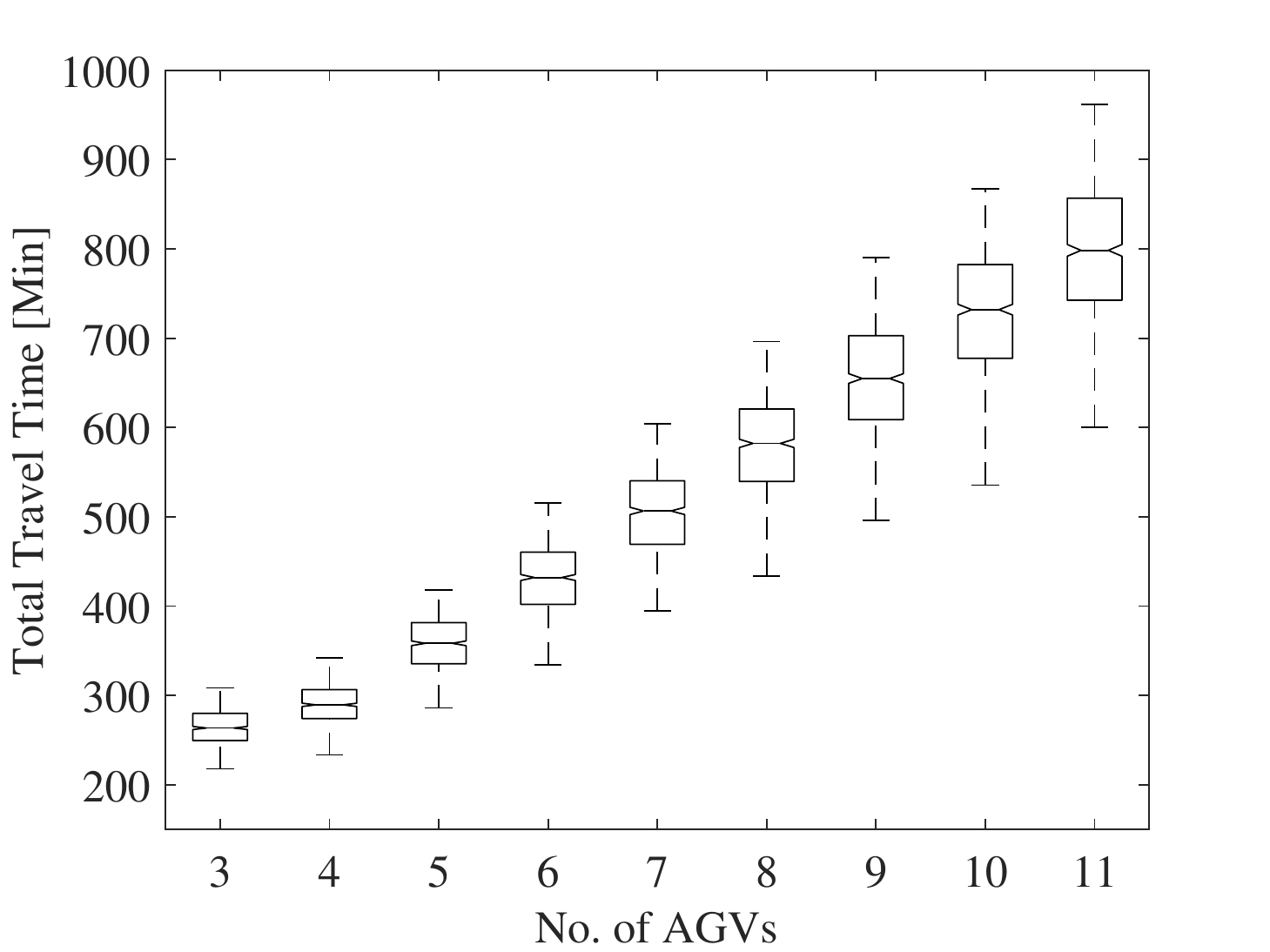}
       \caption{Total Travel Time Model S}
        \label{fig:TTTCM2}
    \end{subfigure}
    \label{fig:CleanCartTravelTimeComparison}
\end{figure}

Figure \ref{fig:CleanCartTravelTimeComparison} summarizes the results of the sensitivity analysis for clean case carts. Figures \ref{fig:TTTCM1} and \ref{fig:TTTCM2} show the box plots of the total daily travel times for models M and S, respectively. Observations show that, in model M, the total travel time increases with the number of AGVs. Since every AGV travels exactly the same distance, the increase in travel time is due to waiting in traffic. Traffic congestion increases with the number of AGVs in the system. Similarly, the travel time is sensitive to the number of AGVs used in model S (see Figure  \ref{fig:TTTCM2}).

A statistical analysis was conducted to test whether the travel times of models M and S are significantly different at a 95\% confidence level. Table \ref{Table:MvsS} shows the results of the $t$-test when the number of AGVs in both the systems is 11. The difference between the average travel times is statistically significant, and the travel time for clean case carts increases if the elevators are swapped. Longer travel times occur because the travel distances for clean carts in the second floor are longer from the new detent area to the ORs.

\begin{table}[H]
	\centering
	\caption{Clean Case Cart Movement: Model M vs Model S} \label{Table:MvsS}
	{\begin{tabular}{ccc}
			\hline
			\textbf{Systems} & \textbf{Average Travel Time} & \textbf{Confidence Interval} \\ \hline
			Model M & 9.67& (9.65,9.68)\\
			Model S & 10.72& (10.72,10.73) \\ \hline
	\end{tabular}}
\end{table}

\begin{figure}[H]
	\centering
	\caption{Sensitivity Analysis Task Completion Times: Model M vs Model S}
	\begin{subfigure}[b]{0.45\textwidth}
		\includegraphics[width=\textwidth]{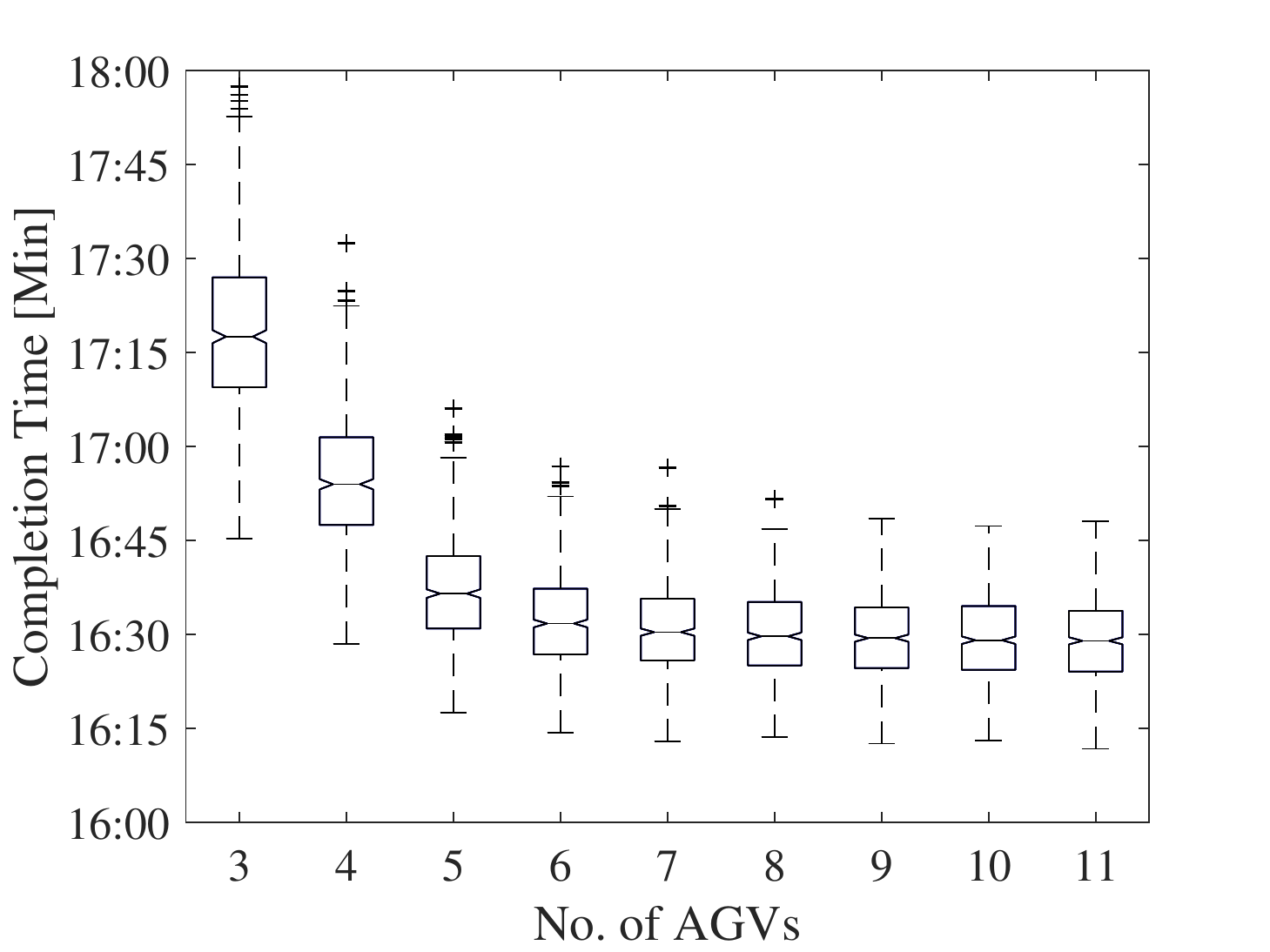}
		\caption{Task Completion Time Model M}
		\label{fig:ACTM1}
	\end{subfigure}
	~ 
	\begin{subfigure}[b]{0.45\textwidth}
		\includegraphics[width=\textwidth]{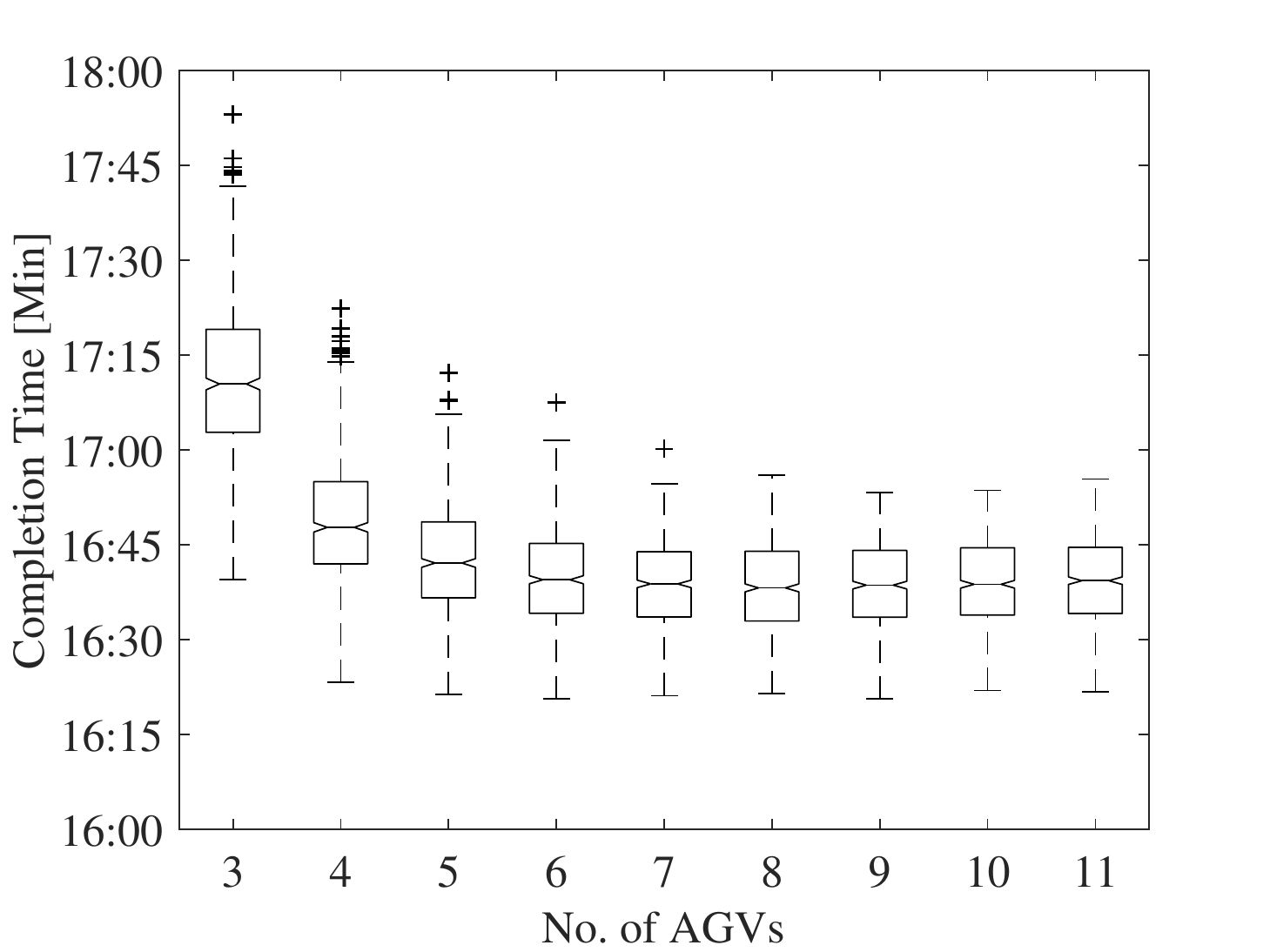}
		\caption{Task Completion Time Model S}
		\label{fig:ACTM2}
	\end{subfigure}
	~ 
	\label{fig:CompletionTimeComparison}
\end{figure}

\textcolor{black}{Figures \ref{fig:ACTM1} and \ref{fig:ACTM2} present the task completion time for models M and S. Task completion time decreases as the number of AGVs increase. A significant reduction in total travel times can be observed in both systems if 6 or 7 AGVs are used. However, increasing the number of AGVs beyond 6 or 7 does not have a significant impact on task completion time.}

\subsubsection{Results of the Sensitivity Analysis: Soiled Cart Movement}

\begin{figure}[H]
    \centering
    \caption{Sensitivity Analysis Soiled of Case Carts: Model M vs Model S}
    \begin{subfigure}[b]{0.45\textwidth}
        \includegraphics[width=\textwidth]{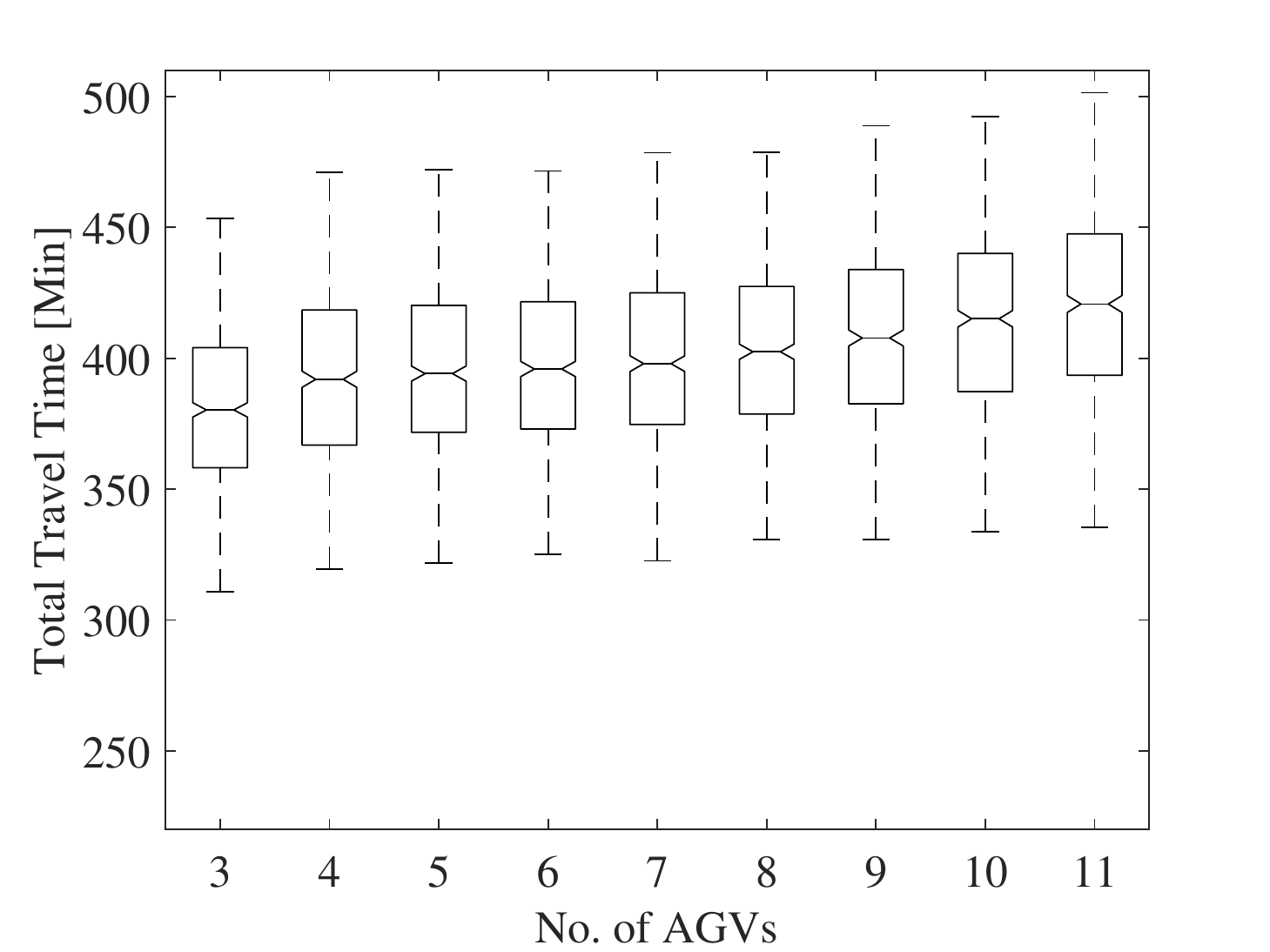}
        \caption{Total Travel Time Model M}
        \label{fig:TTTDM1}
    \end{subfigure}
    ~ 
    \begin{subfigure}[b]{0.45\textwidth}
        \includegraphics[width=\textwidth]{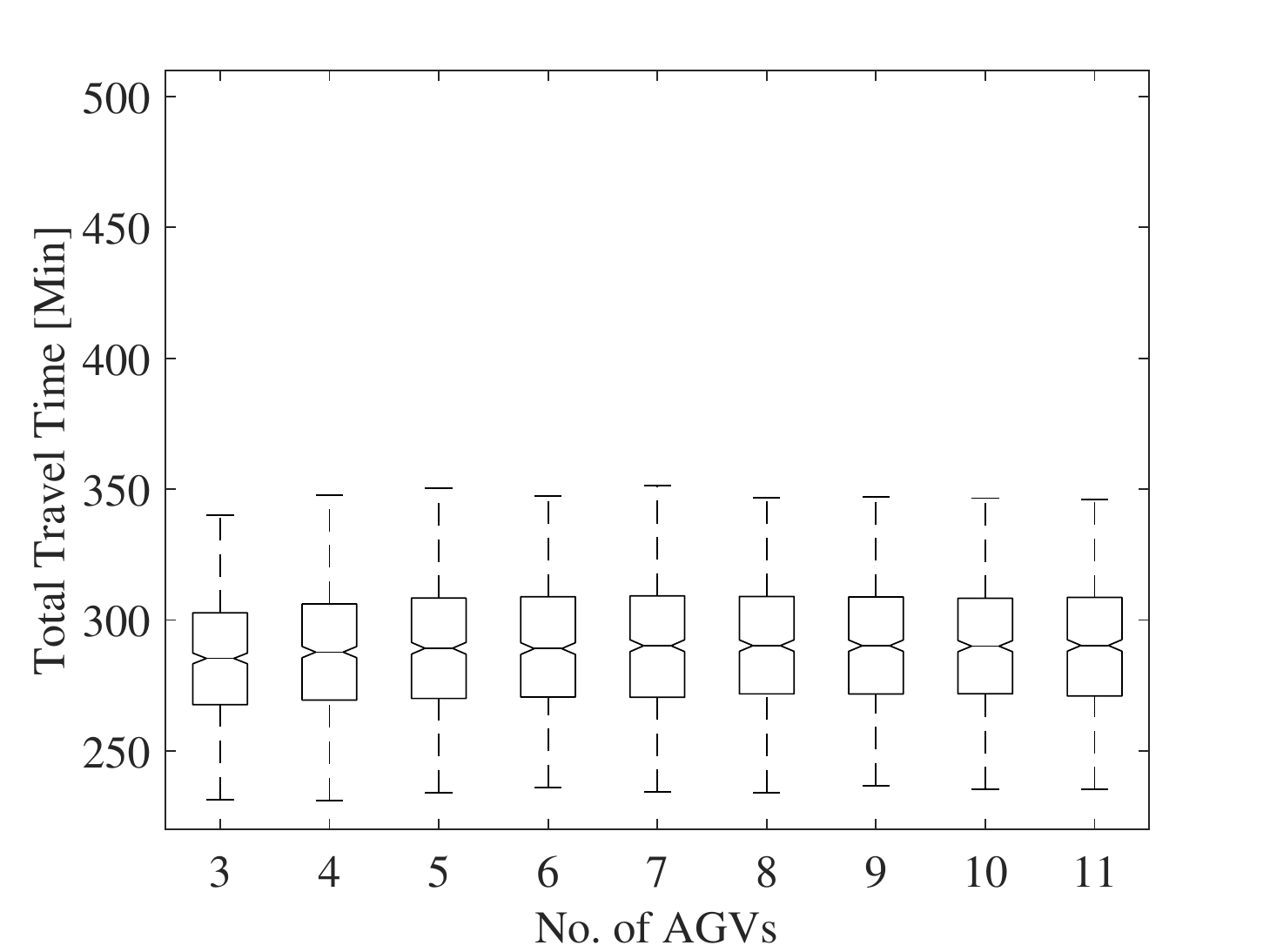}
        \caption{Total Travel Time Model S}
        \label{fig:TTTDM2}
    \end{subfigure}
    ~ 
    \label{fig:SoiledCartTravelTimeComparison}
\end{figure}

Figure \ref{fig:SoiledCartTravelTimeComparison} summarizes the results of the sensitivity analysis for soiled case carts. Figures \ref{fig:TTTDM1} and \ref{fig:TTTDM2} present box plots of the total daily travel times for models M and S respectively. In model M, an increase in the number of AGVs leads to longer travel times \textcolor{black}{for soiled case carts} because of congestion.  This is mainly because AGVs with clean and soiled carts share paths.  Travel times of soiled carts in model S are not impacted by changes to the  number of AGVs. 


It was tested whether the travel times in model M are significantly different than the travel times in model S for soiled case carts, using a two-sample t-test at a 95\% confidence level. Table \ref{Table:MvsSSoiled} shows the results of the $t$-test when the number of AGVs in both the systems is 11. It is evident that the difference between the average travel times is statistically significant, and the travel time for soiled case carts decreases if the elevators are swapped.

\begin{table}[H]
\centering
\caption{Soiled Case Cart Movement: Model M vs Model S}
\label{Table:MvsSSoiled}
{\begin{tabular}{ccc}\hline
\textbf{Systems} & \textbf{Average Travel Time} &\textbf{ Confidence Interval} \\\hline
Model M & 6.23& (6.22,6.24) \\
Model S & 4.26& (4.25,4.26) \\ \hline
\end{tabular}}
\end{table}


\subsection{Research Question 2: Improve the Process Via a Kanban System with Limited Number of AGVs}

The data analysis shows that the volume of surgical cases follows a distribution whose mean value changes based on the day of the week. Hence, one would expect that the hospital could use a different number of AGVs to deliver surgical case carts on each day, \textcolor{black}{instead of using a fixed number of AGVs every day}. Such an approach could lead to improved AGV utilization, reduced congestion, and shorter trip times.

A practical solution is proposed, based on the principles of a Kanban system, to control the number of AGVs used daily. In particular, consider a k-container Kanban system with two processes: the loading of a new surgical case cart and the consumption of the case cart at the destination. The Kanban system operates with the AGV as the container. A card is attached to each AGV, and an empty container is returned to the pickup location for replenishment. Notice that, instead of the card being attached to the materials, the container becomes the actual Kanban, in this case, the AGVs. Also, the container is available or empty after the case cart has been dropped and is returned to the pickup location for the replenishment, i.e., to pick up another case cart.

Two experiments are conducted using ARENA OptQuest. Both experiments focus on optimizing AGV movements on ``\textit{Path of AGVs with Clean Cart"} since analysis shows the movement of clean carts causes congestion on the mezzanine floor. For both experiments, the decision variable is the \textit{number of AGVs to be used on a particular day of the week}.

In experiment 1, the objective is to minimize the total travel time each day. It is also important that the movements of all clean case carts are completed by 7 pm. To ensure that, in experiment 1, a task completion time constraint $T_c \leq 200$ minutes is added. In experiment 2, the objective is to minimize the sum of task completion times over a replication. The results of the OptQuest output are discussed in the next section.

Based on the results of simulation-optimization experiments for a given day of the week, the number of AGVs to be used for that day is determined. The practical implementation of this system is as follows. Suppose there are $N$ AGVs in the system dedicated to surgical case cart movement, and simulation optimization recommends the use of $k\leq N$ AGVs on a particular day. If all $k$ AGVs are active, i.e., all Kanban containers are in use, the system prevents requesting additional AGVs. When an AGV is freed at a drop-off location, it is made available to the system again. The Kanban system provides GMH with an easy way of limiting and monitoring the number of AGVs used each day. The implementation does not require any additional cost except for minor changes in the control system. The Kanban system can also be implemented by simply limiting the number of AGVs requested. A new AGV may not be requested until one of the active AGVs has completed its task. This method requires a high level of communication between departments.

\subsubsection{Results of the Kanban System Implementation in OptQuest}\label{OptQuestScenarios}

\textcolor{black}{Tables \ref{Table:Obj1} and \ref{Table:Obj2} summarize the results of the OptQuest experiments. These tables present the solutions that satisfy the following three conditions: \textit{(i)} the total number of AGVs used in a day is less than or equal to 11; \textit{(ii)} the average travel time per AGV is less than or equal to the average travel time observed from the data; and \textit{(iii)} the total completion time is no later than 5:05 pm. These criteria identify solutions that could potentially be adopted by GMH.  Each solution presents the minimum, maximum, and average travel time for each AGV; the task completion time; and the number of AGVs used each day.}


The results in Table \ref{Table:Obj1} suggest the use of fewer AGVs than the current practice at GMH because the objective of experiment 1 is to minimize the total travel time. Utilizing fewer AGVs leads to reduced congestion, as evidenced by the average travel time and the corresponding range of travel time, which is narrower. In contrast, when the objective is to minimize the task completion time, the simulation experiments suggest using relatively more AGVs, \textcolor{black}{as can be seen from the results in Table \ref{Table:Obj2}}. This increase in the number of AGVs leads to congestion, evidenced by the average travel time and the corresponding range of travel time, which is wider compared to results in Table \ref{Table:Obj1}. Note that completion time is impacted by travel time and waiting time. As a result, the task completion time is \textcolor{black}{shorter for solutions with a higher number of AGVs available} because carts do not wait for another AGV.

\textcolor{black}{The solutions of OptQuest use a different number of AGVs each day of the week, which is different from the current practice at GMH. Experimental results suggest that, on days with a lower case volume, fewer AGVs should be used than on days with a higher volume of cases.}

\begin{table}[H]
\centering
\caption{Results of Experiment 1: Minimize Total Travel Time Per Day}
\label{Table:Obj1}
{\footnotesize \begin{tabular}{c|ccc|c|ccccc}
\cline{1-10}
 & \multicolumn{3}{c|}{\textbf{Travel Times {[}Min{]}}} & \textbf{Task Completion Time } & \multicolumn{5}{c}{\textbf{Number of AGVs}} \\ \cline{1-10} 
 \textbf{\textcolor{black}{Solution}} & \textbf{Min} & \textbf{Max} & \textbf{Average} & \textbf{Average} & \textbf{M} & \textbf{T} & \textbf{W} & \textbf{Th} & \textbf{F} \\ \hline
 1 & 2.64 & 6.47 & 3.43 & 5:03:01 PM & 3 & 3 & 3 & 4 & 4 \\ 
 2 & 2.64 & 6.78 & 3.58 & 4:55:51 PM & 3 & 4 & 4 & 3 & 5 \\ 
 3 & 2.64 & 9.72 & 5.37 & 4:41:20 PM & 3 & 4 & 8 & 7 & 8 \\ 
 4 & 2.64 & 12.1 & 5.35 & 4:37:50 PM & 4 & 4 & 5 & 10 & 7 \\ 
 5 & 2.64 & 6.78 & 3.84 & 4:43:25 PM & 5 & 4 & 4 & 4 & 5 \\ 
 6 & 2.64 & 7.93 & 4.33 & 4:34:51 PM & 5 & 4 & 6 & 5 & 6 \\ \hline
\end{tabular}}
{}
\end{table}

\begin{table}[H]
\centering
\caption{Results of Experiment 2: Minimize the Task Completion}
\label{Table:Obj2}
{\footnotesize \begin{tabular}{c|ccc|c|ccccc}
\cline{1-10}
 & \multicolumn{3}{c|}{\textbf{Travel Times {[}Min{]}}} & \textbf{Task Completion Time} & \multicolumn{5}{c}{\textbf{Number of AGVs}} \\ \cline{1-10} 
 \textbf{\textcolor{black}{Solution}} & \textbf{Min} & \textbf{Max} & \textbf{Average} & \textbf{Average} & \textbf{M} & \textbf{T} & \textbf{W} & \textbf{Th} & \textbf{F} \\ \hline
 1 & 2.64 & 12.1 & 7.27 & 4:28:46 PM & 8 & 7 & 8 & 10 & 8 \\ 
 2 & 2.64 & 12.1 & 7.96 & 4:28:30 PM & 10 & 7 & 10 & 10 & 8 \\ 
 3 & 2.64 & 13.64 & 8.68 & 4:28:01 PM & 11 & 6 & 11 & 11 & 11 \\
 4 & 2.64 & 13.64 & 8.57 & 4:28:04 PM & 11 & 7 & 10 & 10 & 11 \\ 
 5 & 2.64 & 13.64 & 8.9 & 4:27:49 PM & 11 & 7 & 11 & 11 & 11 \\ \hline
\end{tabular}}
{}
\end{table}

\section{Implementation}\label{results}

\textcolor{black}{To further evaluate the impact of the proposed Kanban system on AGV utilization, travel time, task completion time and congestion, the solutions obtained from the simulation experiments were implemented using the following approaches: First, a short pilot study was conducted at GMH. This study was only one week long because of the additional resources needed for implementation. Section \ref{PilotStudy} summarizes the results of this study. Next, a second study was conducted via simulation usinf real-life data from GMH regarding the total number of surgical cases conducted each day of the week, from January 1, 2018, through September 11, 2018. This 26 weeks worth of data allowed a thorough statistical analysis of the results. A fleet of AGVs was selected to deliver the surgical carts each week of this period. Section \ref{SimImplement} summarizes the results of this study. Section \ref{ManagerialInsights} presents the managerial insights revealed by the pilot study and the corresponding simulation experiments.}

\subsection{A Pilot Study at GMH}\label{PilotStudy}

\textcolor{black}{Our proposed Kanban system was piloted at GMH for one week. During the pilot study, we visited GMH every day and collected data on the movement of AGVs from 3:45 pm to 5 pm. Throughout this week, the number of AGVs used at GMH every day was the same as the first solution presented in Table  \ref{Table:Obj2}. A less conservative solution with more moderate number of AGVs was used instead of the solutions presented in Table \ref{Table:Obj1} because GMH staff expressed concerns about the potential delays that may result from significantly reducing the number of AGVs used each day (from 11 to 3, 4, or 5). Of all the solutions presented in Table \ref{Table:Obj2}, solution 1 uses the fewest number of AGVs in a day. Therefore, this implementation allowed us to evaluate how reducing the number of AGVs would impact congestion.}

\begin{table}[H]
\centering
\caption{Pilot Study Results: Average Travel Times} \label{Table:PilotStudy2}
{\footnotesize \begin{tabular}{l|cc|cc|cc}
\hline
&\multicolumn{2}{c|}{\textbf{Week Before}} & \multicolumn{2}{c|}{\textbf{\textcolor{black}{Treatment Week}}} & \multicolumn{2}{c}{\textbf{Week After}} \\ 
\cline{2-7}
\multicolumn{1}{c|}{\textbf{Day}}  &  & \textbf{Avg. } &  & \textbf{Avg.} &   & \textbf{Avg.} \\ 
 &\multicolumn{1}{c}{\textbf{Case Vol.}} & \textbf{Travel Time} & \textbf{Case Vol.} & \textbf{Travel Time} & \textbf{Case Vol.} & \textbf{Travel Time} \\ \hline
Monday & \multicolumn{1}{c}{28} & 9.07 & 26 & 5.50 & 31 & 11.37 \\ 
Tuesday & \multicolumn{1}{c}{23} & 17.13 & 34 & 9.53 & 21 & 16.71 \\ 
Wednesday & \multicolumn{1}{c}{30} & 8.57 & 30 & 7.73 & \textbf{14} & 5.35 \\ 
Thursday & \multicolumn{1}{c}{24} & 16.79 & 26 & 8.54 & 25 & 9.35 \\ 
Friday & \multicolumn{1}{c}{30} & 10.37 & 32 & 8.00 & 22 & 9.26 \\ \hline
\end{tabular}}
{}
\end{table}

Table \ref{Table:PilotStudy2} presents the results of actual travel times during \textcolor{black}{3:45 pm to 5 pm each day of the treatment} week ($t$), the week before ($t-1$), and week after ($t + 1$). The travel times during weeks $t-1$ and $t+1$ are tested, using a two-sample $t$-test at a 95\% confidence level, to determine whether they are significantly different than the travel times during week $t$ for clean case carts. \textcolor{black}{Table \ref{Table:PilotStudy} summarizes the corresponding results of the two-sample $t$-test}.

\begin{table}[H]
\centering
\caption{Pilot Study Results: P-Values}
\label{Table:PilotStudy}
{\begin{tabular}{ccc}
\hline
\multicolumn{1}{c|}{\textbf{\textcolor{black}{Treatment Period}}} & \multicolumn{2}{c}{\textbf{P-Values (Avg. Travel Times)}} \\ \hline
\multicolumn{1}{c|}{\textbf{Day}}            & \textbf{Week Before}            & \textbf{Week After}             \\ \hline
\multicolumn{1}{c|}{Monday}         & 0.01                & 0.00                   \\ 
\multicolumn{1}{c|}{Tuesday}        & 0.00                    & 0.00                   \\ 
\multicolumn{1}{c|}{Wednesday}      & \textbf{0.50}                    & 0.00                   \\ 
\multicolumn{1}{c|}{Thursday}       & 0.00                    & \textbf{0.60}                   \\ 
\multicolumn{1}{c|}{Friday}         & 0.02                    & \textbf{0.43}\\
\hline               
\end{tabular}}
{}
\end{table}

The average travel time during week $t$ was lower than the average travel time during week $t-1$ on all five days. This difference was statistically significant on Monday, Tuesday, Thursday, and Friday. On the other hand, the average travel time during week $t$ was lower than the average travel time during week $t + 1$ on four days. This difference was statistically significant on Monday and Tuesday. The results show that the average travel time during week $t$ was greater than the week $t + 1$ only on Wednesday, and the difference was statistically significant. This difference can be attributed to the fact that the number of cases on Wednesday in week $t + 1$ was smaller than half the number of cases in the treatment week on the corresponding day. \textcolor{black}{For very low case volumes, the material handling system will not be considerably affected by congestion, and using more AGVs does not have as many adverse effects.}

\begin{table}[h]
\centering
\caption{Pilot Study Results: Standard Deviation of Travel Times}
\label{Table:PilotStudy3}
{\begin{tabular}{cccc}
\hline
\textbf{\textcolor{black}{Treatment Period}} & \multicolumn{3}{c}{\textbf{Std. Deviation }} \\ \hline
\textbf{Day} & \textbf{Week Before} & \textbf{\textcolor{black}{Treatment Week}} & \textbf{Week After} \\ \hline
Monday & 6.97 & 1.72 & 3.82 \\
Tuesday & 9.67 & 3.34 & 9.02 \\
Wednesday & 6.06 & 3.46 & \textbf{1.34} \\
Thursday & 7.53 & 2.58 & 7.12 \\
Friday & 4.83 & 2.68 & 7.19 \\ \hline
\end{tabular}}
{}
\end{table}

Table \ref{Table:PilotStudy3} presents the standard deviation of travel times during the weeks $t-1$, $t$, and $t+1$. The standard deviation of travel times during weeks $t-1$ and $t+1$ is tested at a 95\% confidence level to determine whether they are significantly different than the standard deviation of travel times during week $t$ for clean case carts. Table \ref{Table:PilotStudy4} provide the p-values of the tests performed.

\begin{table}[H]
\centering
\caption{Pilot Study Results: P-Values}
\label{Table:PilotStudy4}
{\begin{tabular}{ccc}
\hline
\multicolumn{1}{c|}{\textbf{Treatment Period}} & \multicolumn{2}{c}{\textbf{P-Values (Std. Deviation of Travel Times)}} \\ \hline
\multicolumn{1}{c|}{\textbf{Day}}            & \textbf{Week Before}            & \textbf{Week After}             \\ \hline
\multicolumn{1}{c|}{Monday}         & 0.00                & 0.00                   \\ 
\multicolumn{1}{c|}{Tuesday}        & 0.00                    & 0.00                   \\ 
\multicolumn{1}{c|}{Wednesday}      & 0.04                    & \textbf{0.00}                   \\ 
\multicolumn{1}{c|}{Thursday}       & 0.00                    & \textbf{0.14}                   \\ 
\multicolumn{1}{c|}{Friday}         & 0.00                    & \textbf{0.18}\\
\hline               
\end{tabular}}
{}
\end{table}

The standard deviation of travel times during week $t$ is less than the standard deviation of travel times during the week $t-1$  on all five days. This difference is statistically significant for all five days. The standard deviation of travel time was lower during week $t$ than $t + 1$ on Monday, Tuesday, Thursday, and Friday. This difference was statistically significant on Monday and Tuesday. Similar to the observations related to the average travel time, the standard deviation of travel time on Wednesday in week $t + 1$ was significantly lower than for week $t$.

It is already established that longer travel time indicates longer wait times due to congestion. Similarly, the standard deviation of travel time is a measure of congestion in the system, i.e., a higher standard deviation, while the traveled distance is the same, indicates longer wait times due to congestion. This analysis of the pilot study results clearly shows that congestion was reduced by limiting the number of AGVs in the system, which led to reduced wait times and, consequently, to reduced travel times.

\textcolor{black}{\textbf{Limitations of pilot study:} 
The data collected via pilot study is not extensive due to short period of implementation. In addition to that, during the treatment week, the movement of surgical carts began at about 3:30 to 3:45 on 2 days, so, only about 60\% of the carts were delivered by 5 pm. At 5 pm, the AGVs were assigned to other tasks (e.g., delivery of dinner), so, the remainder of the carts was delivered later on in the evening at about 9 pm, when AGVs were available. Because of this lack of data,  the completion time times are not reported. To overcome these limitations, additional experiments are conducted to test the proposed Kanban system using real-life data for 26 weeks. The implementation via simulation is discussed in the next subsection.}

\subsection{Implementation via Simulation}\label{SimImplement}

\begin{figure}[H]
	\centering
	\caption{Evaluation of Policies from Experiment 1: Minimize Total Travel Time Per Day}
	\begin{subfigure}[b]{0.45\textwidth}
		\includegraphics[width=\textwidth]{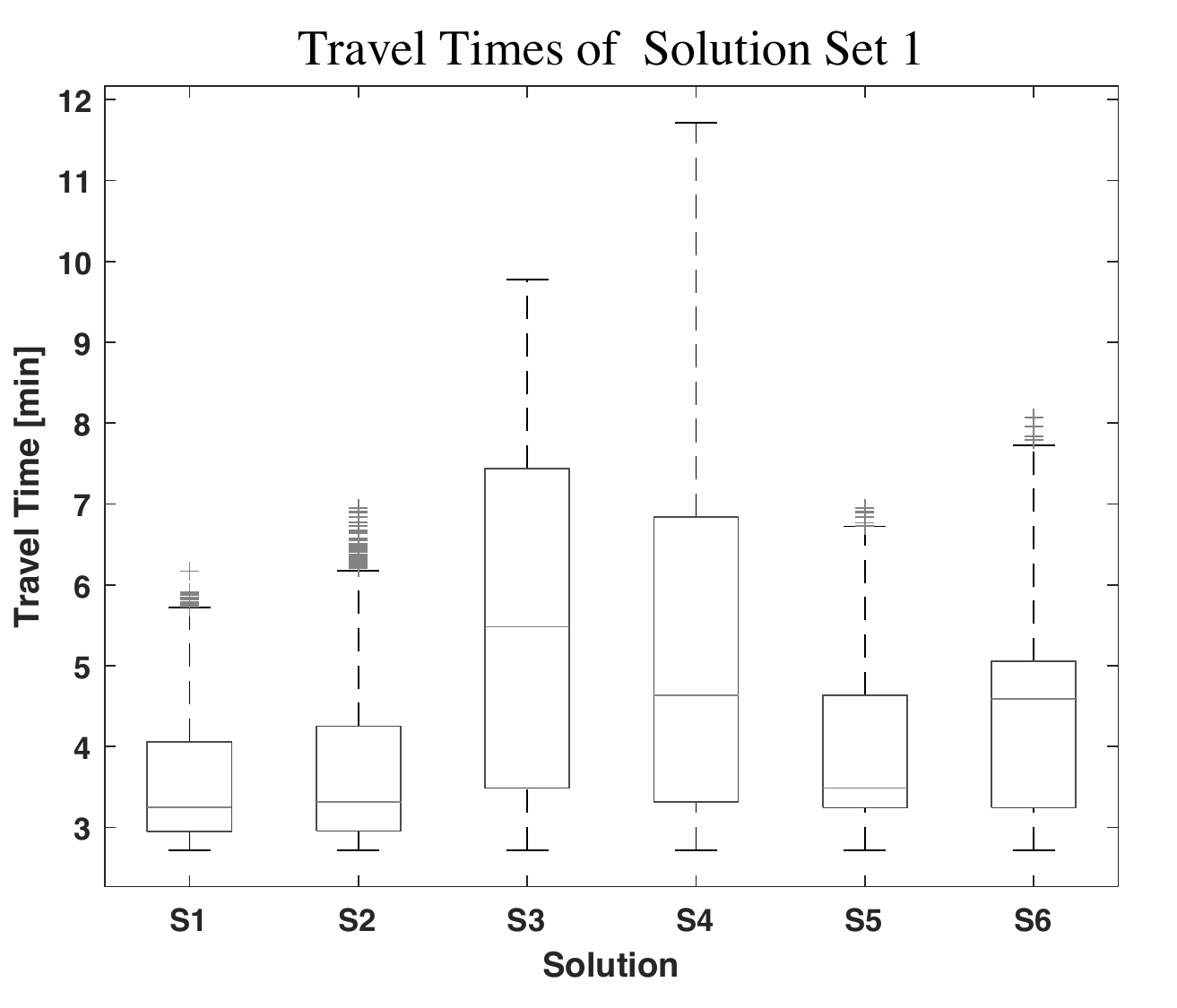}
		\caption{Travel Time}
		\label{fig:Ex1TravelTime}
	\end{subfigure}
	\begin{subfigure}[b]{0.5\textwidth}
		\includegraphics[width=\textwidth]{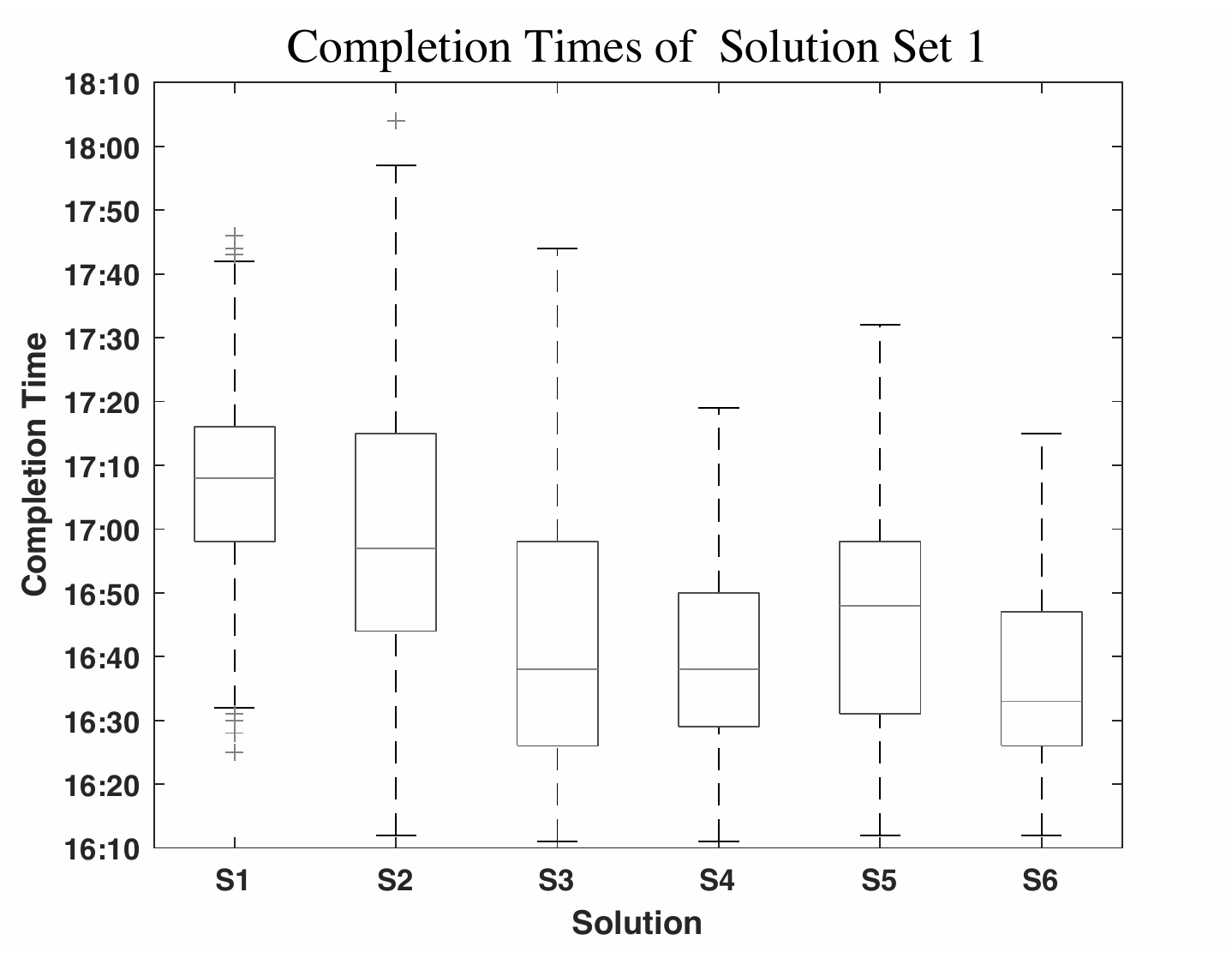}
		\caption{Task Completion Time}
		\label{fig:Ex1CompletionTime}
	\end{subfigure}
	\label{fig:Ex1Evaluation}
\end{figure}

\begin{figure}[H]
	\centering
	\caption{Evaluation of Policies from Experiment 2: Minimize Task Completion Time}
	\begin{subfigure}[b]{0.45\textwidth}
		\includegraphics[width=\textwidth]{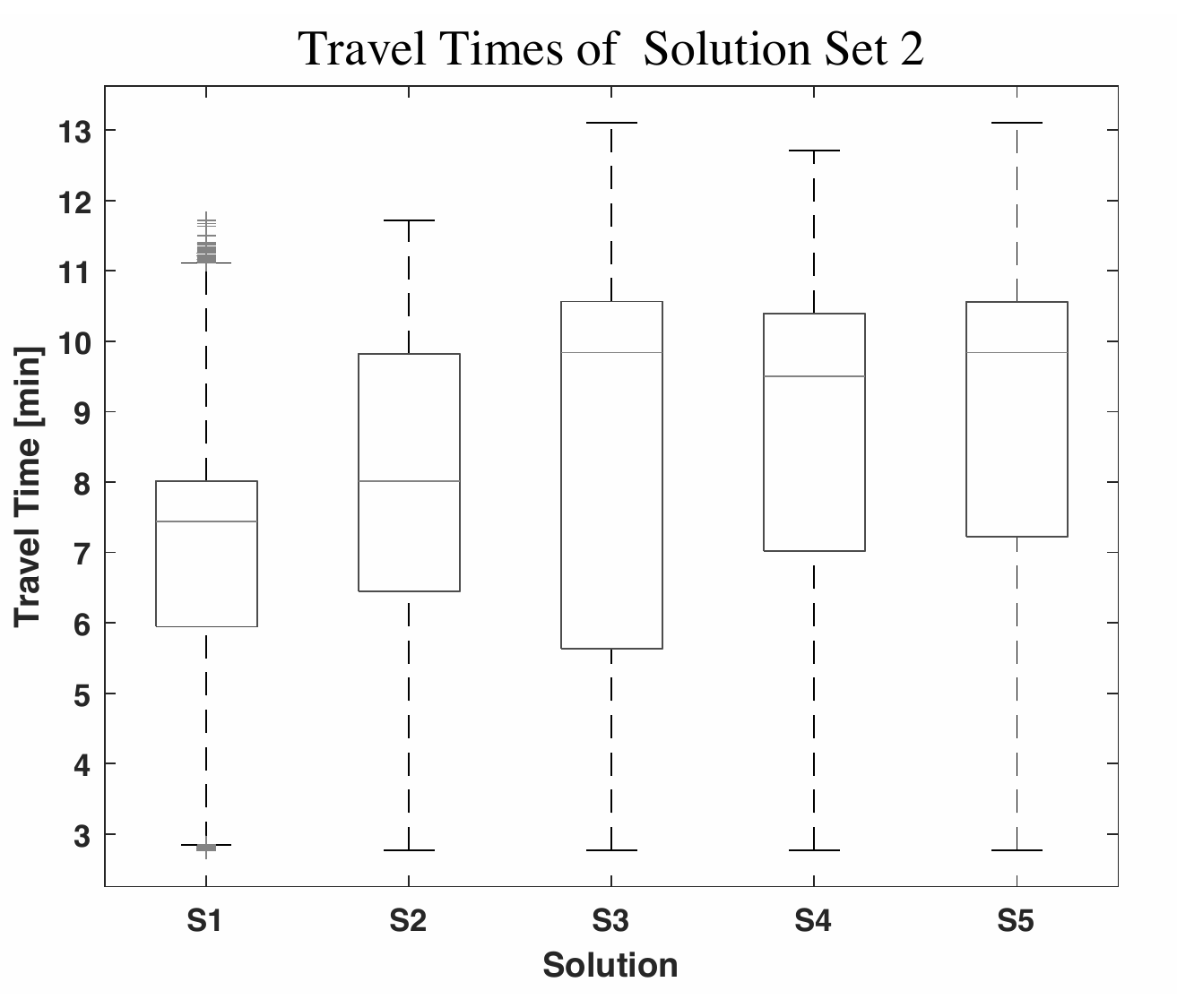}
		\caption{Travel Time}
		\label{fig:Ex2TravelTime}
	\end{subfigure}
	\begin{subfigure}[b]{0.5\textwidth}
		\includegraphics[width=\textwidth]{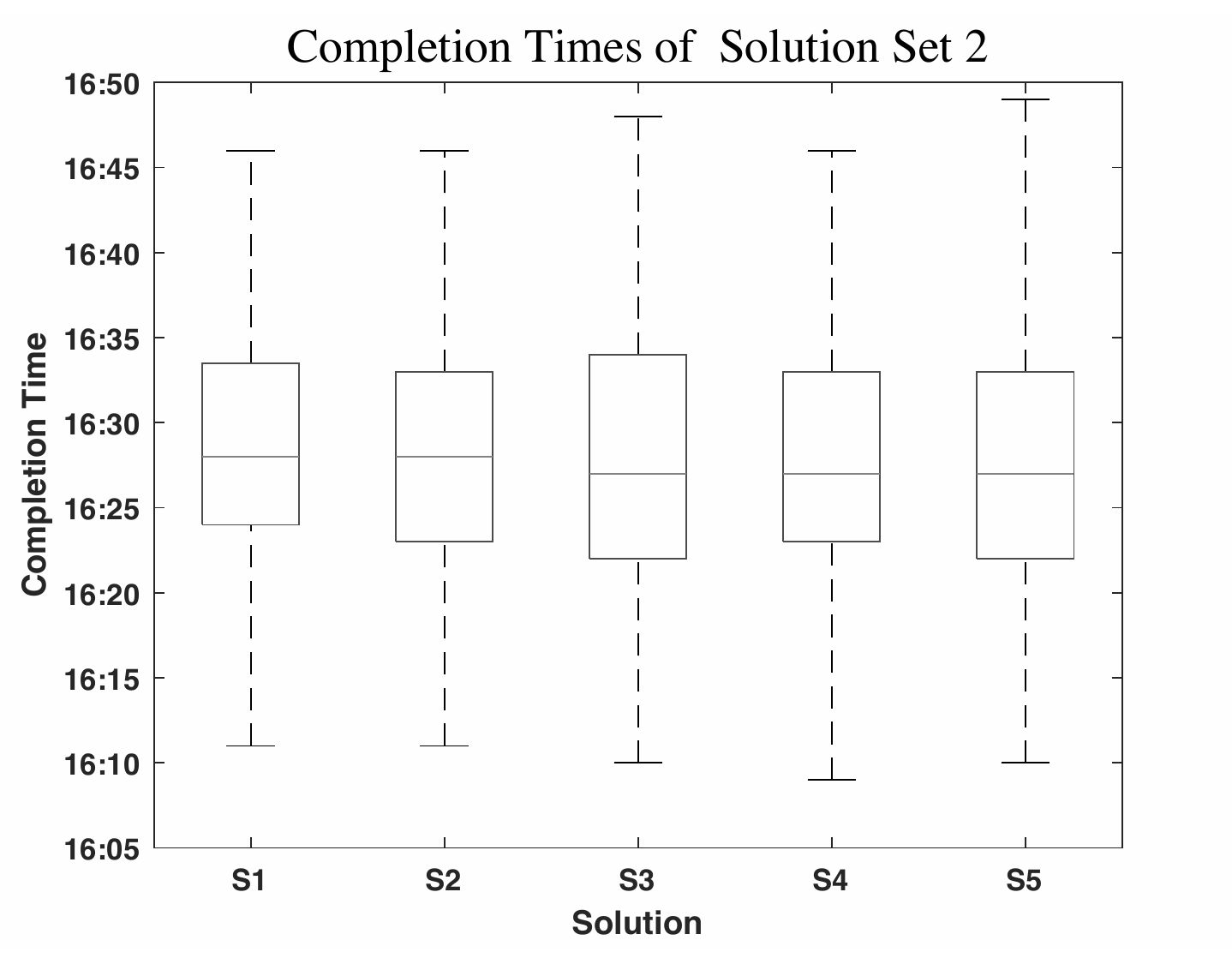}
		\caption{Task Completion Time}
		\label{fig:Ex2CompletionTime}
	\end{subfigure}
	\label{fig:Ex2Evaluation}
\end{figure}

\textcolor{black}{The results of this pilot study indicate that the implementation of a Kanban system where the number AGVs in rotation is limited and varied based on the day of the week has potential to improve travel times. However, the data collected via the pilot study is not extensive due to the short implementation period. To overcome this limitation and also evaluate other solutions generated in Section \ref{OptQuestScenarios}, a set of simulation experiments were run. Actual case volume data, collected by GMH for 26 weeks between January 1, 2018, and September 11, 2018, was used. Figures \ref{fig:Ex1Evaluation} and \ref{fig:Ex2Evaluation} depict the results of the simulation runs for the scenarios obtained from experiment 1 (for solution set 1, see Table \ref{Table:Obj1}) and experiment 2 (for solution set 2, see Table \ref{Table:Obj2}), respectively. Tables \ref{Table:SimImp1} and \ref{Table:SimImp2} present the results of the simulations run to compare the solution implemented in the pilot study with the current practice at GMH.}

\textcolor{black}{The simulation results confirm earlier observations with respect to the simulation -optimization results. Evaluation of the scenarios generated to minimize total travel times shows that the solutions that deploy more AGVs result in higher average travel times but relatively early task completion times. Additionally, all solutions but one result in average task completion times earlier than 5 pm. Comparison of the results depicted in Figure \ref{fig:Ex1Evaluation} with the results reported in Table \ref{Table:Obj1} validate the relationship identified here among the number of AGVs used, congestion levels, and the performance measures considered.}

\textcolor{black}{Similarly, evaluation of the scenarios generated to minimize task completion time confirms the observations presented above. The simulation results show that the solutions considered here yield comparable task completion times which are, on average, around 4:30 pm., well before the target completion time determined by the hospital. On the other hand, the average travel time attained varies among the solutions considered. The robustness of the task completion times can be explained by the tradeoff between travel times and the number of simultaneous trips possible (i.e., the number of AGVs), which is, again, consistent with earlier results.}

\begin{table}[H]
	\centering
	\caption{Implementation via Simulation: Average Travel Times} \label{Table:SimImp1}
	{\footnotesize \begin{tabular}{l|ccc|ccc}
			\hline
			&\multicolumn{3}{c|}{\textbf{Proposed Solution}} & \multicolumn{3}{c}{\textbf{\textcolor{black}{Current Practice}}} \\ 
			\cline{2-7}
			\textbf{Day} & \textbf{Avg. } & &  & \textbf{Avg. }  &  &   \\ 
			&{\textbf{Travel Time}} & \textbf{St.Dev} & \textbf{Confidence Int.} & \textbf{Travel Time} & \textbf{St. Dev} & \textbf{Confidence Int.} \\ \hline
			Monday & 6.82 & 0.036 & (6.810, 6.838)
			 & 9.61 & 0.057 & (9.593, 9.636) \\ 
			Tuesday & 5.84 & 0.036 & (5.827, 5.853) & 9.61 & 0.072 & (9.585, 9.639) \\ 
			Wednesday & 6.45 & 0.036 & (6.434, 6.460) & 8.96 & 0.048 & (8.943, 8.979) \\ 
			Thursday & 7.96 & 0.057 & (7.936, 7.978) & 8.77 & 0.073 & (8.738, 8.792) \\ 
			Friday & 6.55 & 0.041 & (6.533, 6.563) & 9.22 & 0.059 & (9.198, 9.242) \\ \hline
	\end{tabular}}
\end{table}

\begin{table}[H]
	\centering
	\caption{Implementation via Simulation: Average Completion Times} \label{Table:SimImp2}
	{\footnotesize \begin{tabular}{l|ccc|ccc}
			\hline
			&\multicolumn{3}{c|}{\textbf{Proposed Solution}} & \multicolumn{3}{c}{\textbf{\textcolor{black}{Current Practice}}} \\ 
			\cline{2-7}
			\textbf{Day} & \textbf{Avg. } & &  & \textbf{Avg. }  &  &   \\ 
			&{\textbf{Comp. Time}} & \textbf{St.Dev} & \textbf{Confidence Int.} & \textbf{Comp. Time} & \textbf{St. Dev} & \textbf{Confidence Int.} \\ \hline
			Monday & 119.37 & 2.436 & (118.5, 120.3) & 117.91 & 2.483 & (116.98, 118.84) \\ 
			Tuesday & 107.54 & 1.829 & (106.9, 108.2) & 104.96 & 1.725 & (104.32, 105.61) \\ 
			Wednesday & 176.31 & 7.175 & (173.6, 179.0) & 175.18 & 7.302 & (172.45, 177.91) \\ 
			Thursday & 177.73 & 6.021 & (175.5, 180.0) & 177.41 & 6.019 & (175.17, 179.66) \\ 
			Friday & 119.94 & 2.259 & (119.1, 120.8) & 118.52 & 2.281 & (117.67, 119.37) \\ \hline
	\end{tabular}}
\end{table}

\textcolor{black}{To supplement the results obtained from the pilot study, simulations with longer replication length were completed for the solution implemented during the pilot study at GMH. The same 26 weeks of real-life case volume data were used to run these experiments. The simulation results summarized in Table \ref{Table:SimImp1} show that the average travel time under the proposed Kanban system significantly reduces the average travel times on all days of the week. Thus, there is benefit in \textit{(i)} varying the number of AGVs based on the day of the week (i.e., based on the case volume) and \textit{(ii)} using fewer AGVs overall. The proposed solution results in longer average completion times for each day of the week than current practice, but the differences in average completion times are relatively small and are only significant on Monday, Tuesday, and Friday. Furthermore, average completion time under the proposed solution is within the 120-minute period preferred by GMH on all days of the week.}

\subsection{Managerial Insights}\label{ManagerialInsights} 
\textcolor{black}{\textbf{Implications of Our Findings:} This research was motivated by inefficiencies in the material handling system at GMH. During the afternoon hours of 3 to 5, AGVs with clean and dirty case carts used this corridor which led to increased congestion and longer trip times. GMH staff was interested in developing analytical solutions that would lead to reduced congestion in the main corridor of the mezzanine floor.}

\textcolor{black}{\textbf{Research Question \textit{(i)}}: The GMH staff suspected that a change of roles from elevator J to elevators G and K would reduce congestion on the mezzanine floor. The research team conducted an extensive data analysis of trips made by AGVs to develop a simulation model in which the role of elevator J was swapped with G and K. The results of the simulation model indicated lower congestion in the mezzanine floor, supporting the GMH staff’s intuition. However, the simulation model also indicated longer overall trip times. In the current system, elevators G and K carry AGVs with soiled case carts, as well as trash and soiled linen carts, to the ground floor. If the elevators are swapped, as per safety guidelines, soiled linen and trash carts must be delivered using different elevators. A number of elevators were considered for these movements. However, the use of these alternative elevators increases the distance traveled by the AGVs to reach an appropriate detent area. It also increases the distance traveled by an employee to move the soiled linen and trash carts to the new detent area. In addition, for the trash and soiled linen carts to use elevator J, a new AGV guide-path must be installed in front of the elevator on the ground floor. The installation of guide-paths is expensive. Despite the benefit of reducing congestion in the mezzanine floor, the swapping of elevators was found to be costly and difficult to implement.
Even though the solution was not implemented, the GMH staff found the results of our model beneficial since \textcolor{black}{they revealed tradeoffs and challenges that GMH staff had not foreseen.}}

\textcolor{black}{\textbf{Research Question \textit{(ii)}}: The results of the simulation-optimization model indicate that reducing the number of AGVs used each day and changing the number of AGVs based on the volume of cases would lead to reduced congestion, shorter trip times, and shorter task completion time. The idea of reducing the number of AGVs used daily was received with doubt by the staff who loaded the carts, because they were concerned that it would lead to delays in delivery time, and, in consequence, overtime work. The management supported the idea, but the implementation of the Kanban system requires updates of the software that governed the movement of AGVs. These updates are completed by FMC-Technology at a fixed cost.}

\textcolor{black}{\textbf{Limitations of Our Findings:} A limitation of this study is that the model developed here focuses only on the movement of AGVs that deliver surgical case carts. This model can be extended to consider the movement of other AGVs as well. Increasing the scope of the model would result in more accurate modeling of \textcolor{black}{AGV availability, AGV traffic, and interactions among different services that use AGVs.}}

\section{Discussion and Conclusions}

Inefficiencies observed in the material handling system at GMH motivated this research. Hospital staff reported long lines of AGVs waiting for the elevator on the mezzanine floor after clean cart delivery begins in the afternoon. In the current layout, clean and soiled case cart movements share a path. Congestion on this shared path results in delays of delivering clean case carts to the OR. It also prevents soiled case carts and instruments being cleaned in a timely manner. Finally, the AGV traffic on the mezzanine floor poses a serious risk to pedestrians.

This paper proposes a framework that integrates data analysis with system simulation and optimization in order to address the following research questions: \textit{(i)} What are the implications of redesigning a hospital’s material handling system? \textit{(ii)} What are the implications of improving a hospital’s material handling process? Data was collected at GMH, and a thorough data analysis was conducted to understand the inefficiencies of the system. This analysis indicated that the number of trips, the average travel time per trip, and the corresponding standard deviation of travel times are higher during the afternoon, which is the time when clean surgical carts are delivered. This analysis provided the necessary data- and knowledge-base for the development of a discrete event simulation model.

In support of research question \textit{(i)}, a simulation-optimization model was developed that swapped the roles of elevator J with elevators G and K. As GMH staff expected, this swap led to reduced congestion in the mezzanine floor. However, the overall trip time increased because of longer travel times along the second floor where the ORs are located. The swap of the elevators also resulted in longer travel times for the AGVs that deliver dirty linen and trash. Finally, additional investments to install AGV guide-paths on the new routes would be required to implement this redesign.       

In support of research question \textit{(ii)}, a simulation-optimization model was developed in which the number of AGVs used during a particular day changed based on the volume of surgical cases per day. The system used the Kanban principles to control the number of AGVs assigned to the delivery of surgical carts. A sensitivity analysis was conducted to evaluate the impact that reducing the number of AGVs has on travel time and task completion time. This analysis indicates that using fewer AGVs is sufficient to complete the daily delivery of surgical carts. For example, using more than 6 AGVs did not have a significant impact on task completion time. However, the travel time per trip increases with the number of AGVs due to increased congestion.

In order to validate the results of the simulation-optimization model, a one-week long pilot study was conducted, during which the number of AGVs used at GMH was determined based on the proposed model. The travel times of the  treatment week were compared with those of the previous week and the week after, during which the number of AGVs in the system was not controlled. This comparison indicates that travel times were shorter during the week-long pilot study. Since the pilot study had a short duration, a longer implementation was conducted via simulation. The model was run using real-life data from GMH for the period January 1, 2018, to September 11, 2018. The number of AGVs used per day was fixed based on the results of the proposed framework. The statistical analysis of the results indicated that the proposed model led to reduced congestion and shorter trip times.

\textcolor{black}{The extensive data analysis and simulations presented here reinforce what the GMH staff already suspected: that GMH is using more AGVs for the delivery of surgical case carts than it needs to, and the number of AGVs should change daily based on the volume of surgical cases. The GMH staff have considered the models and experimental results as valuable inputs, and have implemented our recommendations in some capacity. These recommendations also played a substantial part in helping GMH evaluate potential changes to be made in the design of the facility. Further improvements can continue to be made in GMH’s material handling systems as the hospital examines other concepts from the manufacturing domain.}

\bibliographystyle{plain}
\bibliography{FinalPaper}

\section*{Appendix}
\begin{table}[H]
	\centering \scalefont{0.88}
\hspace*{-1cm}	\begin{tabular}{p{1cm}|p{16.5cm}}
		\hline
		Day & Distribution \\ \hline
		Mon.  & DISC(0.004,0,0.05,30,0.096,60,0.139,90,0.181,120,0.248,150,0.309,180,0.355,210,0.397,240,0.444,270,\\ \qquad & 0.516,300,0.562,330,0.614,360,0.662,390,0.704,420,0.73,450,0.758,480,0.789,510,0.828,540,0.861,570,\\ \qquad & 0.889,600,0.895,630,0.917,660,0.928,690,0.932,720,0.946,750,0.952,780,0.954,810,0.959,840,0.961,870,\\ \qquad & 0.965,900,0.969,930,0.969,960,0.972,990,0.976,1020,0.976,1050,0.978,1080,0.98,1110,0.983,1140,\\ \qquad & 0.987,1170,0.987,1200,0.989,1230,0.993,1260,0.993,1290,0.998,1320,0.998,1350,0.998,1380,1,1410) \\
		Tue. & DISC(0.006,0,0.035,30,0.088,60,0.146,90,0.203,120,0.259,150,0.307,180,0.359,210,0.42,240,0.461,270,\\ \qquad & 0.518,300,0.572,330,0.61,360,0.653,390,0.689,420,0.72,450,0.752,480,0.789,510,0.821,540,0.85,570,\\ \qquad & 0.875,600,0.902,630,0.908,660,0.919,690,0.929,720,0.935,750,0.946,780,0.948,810,0.952,840,0.958,870,\\ \qquad & 0.96,900,0.965,930,0.969,960,0.969,990,0.973,1020,0.975,1050,0.975,1080,0.975,1110,0.979,1140,\\ \qquad & 0.985,1170,0.987,1200,0.988,1230,0.99,1260,0.992,1290,0.994,1320,0.996,1350,1,1380,1,1410) \\
		Wed. & DISC(0,0,0.004,30,0.033,60,0.075,90,0.133,120,0.18,150,0.264,180,0.308,210,0.353,240,0.399,270,\\ \qquad & 0.472,300,0.523,330,0.577,360,0.621,390,0.645,420,0.694,450,0.74,480,0.785,510,0.818,540,0.843,570,\\ \qquad & 0.893,600,0.914,630,0.923,660,0.934,690,0.944,720,0.951,750,0.958,780,0.964,810,0.964,840,0.964,870,\\ \qquad & 0.969,900,0.974,930,0.976,960,0.978,990,0.978,1020,0.984,1050,0.985,1080,0.985,1110,0.989,1140,\\ \qquad & 0.991,1170,0.995,1200,0.995,1230,0.995,1260,0.995,1290,0.995,1320,0.995,1350,0.998,1380,1,1410) \\
		Thu. & DISC(0.006,0,0.043,30,0.111,60,0.159,90,0.203,120,0.263,150,0.31,180,0.355,210,0.413,240,0.462,270,\\ \qquad & 0.513,300,0.544,330,0.592,360,0.638,390,0.675,420,0.713,450,0.754,480,0.787,510,0.818,540,0.838,570,\\ \qquad & 0.874,600,0.895,630,0.912,660,0.925,690,0.933,720,0.939,750,0.946,780,0.952,810,0.953,840,0.956,870,\\ \qquad & 0.959,900,0.963,930,0.969,960,0.97,990,0.97,1020,0.973,1050,0.976,1080,0.979,1110,0.98,1140,\\ \qquad & 0.98,1170,0.98,1200,0.982,1230,0.984,1260,0.989,1290,0.993,1320,0.993,1350,0.994,1380,1,1410) \\ 
		Fri. & DISC(0.016,0,0.053,30,0.105,60,0.156,90,0.2,120,0.268,150,0.353,180,0.411,210,0.451,240,0.486,270,\\ \qquad & 0.53,300,0.579,330,0.626,360,0.674,390,0.716,420,0.751,450,0.788,480,0.805,510,0.844,540,0.875,570,\\ \qquad & 0.9,600,0.912,630,0.923,660,0.933,690, 0.94,720,0.951,750,0.954,780,0.954,810,0.956,840,0.96,870,\\ \qquad & 0.963,900,0.963,930,0.967,960,0.967,990,0.97,1020,0.97,1050,0.974,1080,0.974,1110,0.975,1140,\\ \qquad & 0.975,1170,0.979,1200,0.981,1230,0.986,1260,0.991,1290,0.993,1320,0.995,1350,0.998,1380,1,1410) \\ \hline
	\end{tabular}\hspace*{-1cm}%
	\caption{Discete Distribution for Soiled Case Cart Release Times}
	\label{tab:discdist}
\end{table}
\end{document}